\input amstex
\input amsppt.sty
\catcode`\@=11
\let\logo@=\relax
\catcode`\@=13
\topmatter
\title
{ T-comitants and the Problem of a Center
 for Quadratic Differential Systems}
\endtitle
\author
    { A.M.Voldman, N.I.Vulpe}
\endauthor
\rightheadtext {the Problem of a Center  for Quadratic  Systems}
\abstract

  The new  necessary and sufficient  affine invariant
conditions for the existence and for determining the number of centers for
general quadratic system are pointed out.
These conditions correspond to the partition of 12-dimensional
coefficient space of indicated system with respect to the number and the
multiplicity of its  finite critical points.

\endabstract
\endtopmatter
\document

  Let us consider  system of differential equations
$$
\eqalign{
   {dx^1\over dt}& =a^{1}+ a^{1}_{\alpha}x^{\alpha} + a^{1}_{\alpha\beta}
    x^{\alpha}x^{\beta}\equiv P_0+P_1+P_2,\cr
   {dx^2\over dt}& =a^{2}+ a^{2}_{\alpha}x^{\alpha} + a^{2}_{\alpha\beta}
    x^{\alpha}x^{\beta}\equiv Q_0+Q_1+Q_2,\
}\eqno(j,\alpha, \beta =1,2)\qquad\quad (1)
$$
where $a^{j}_{\alpha}$ and $a^{j}_{\alpha\beta}\ (j,\alpha, \beta =1,2)$ are
real numbers ( the tensor $a^{j}_{\alpha\beta}$ is symmetric in the lower
indices, with respect to which the complete contraction was made) and
$P_i(x^1,\,x^2)\ (i=0,1,2)$ are homogeneous polynomials of degree $i$.
\smallskip

  The existence of a center at the origin for system (1)
was examined for the first time in 1904 by H.Dulac [1]. More precisely,
in [1] the problem of a center was considered for the
following equation:
$$
   {dy\over dx}=-{x+ax^2+bxy+cy^2\over y+a^\prime x^2+b^\prime xy+
   c^\prime y^2},\eqno(2)
$$
with the complex variables $x$ and $y$ and the complex coefficients
$a,\,b,\,c,\,a^\prime,\,b^\prime,\, c^\prime$.

  In W.Kapteyn's papers [2,3] the analogical problem was examined for the
equation (2)  with real variables and coefficients.

 However, neither Dulac nor Kapteyn have been obtained  explicit conditions,
through the coefficients of equation (2), which ensure the existence of
a center at the origin. This problem, besides the determination
of the qualitative phase portraits of equation (2) with a center, were stated
by M. Frommer in [4]. But, as was shown by N.A.Saharnicov [5], there are some
mistakes in  Frommer's paper.

   K.S.Sibirsky [6] and L.N.Belyustina [7] are the first who found out
the explicit conditions,  expressed
through the coefficients of equation (2)
 for the existence of a center. We remark, that the center problem
for equation (2) was also examined by K.E.Malkin [8], I.S.Kukles [9]
and other authors.

   Some years later,  applying
the developed  theory of algebraic invariants of differential
equations, K.S.Sibirsky  [10]  solved the problem under consideration
for a more general case. The center affine invariant conditions for the
existence of a center at the origin  are determined for system (1),
i.e. for system with $a^1=a^2=0$.

Finally,  the necessary and sufficient  affine invariant
conditions for the existence and and for determining the number
of centers (anywhere in the plane)
for a general system (1) were established in [11,12]. The obtained conditions
are expressed as
equalities or inequalities involving polynomials, with the degrees at most
  24.   The following question was formulated:
\smallskip

{\bf QUESTION} {\rm [11].} {\sl Is there a semialgebraic solution of this
problem with a lower maximal degree for the polynomials ?}
\smallskip

 In this paper, by using  T-comitants, we have answered  this
question affirmatively.
  We have pointed out the new   necessary and
sufficient  affine invariant
conditions for the existence and for determining the number of centers for
system (1), corresponding to the
partition of the coefficient space $R^{12}$ of a non-degenerated system (1)
with respect to the number and the multiplicity of the finite critical
points of this system [13,14].

\bigskip
 \centerline{\bf PRELIMINARIES}
\bigskip

   Let $a\in R^{12}$ be an element of the space of the coefficients of the
   system (1)
and let us consider the group  $Q$ of nondegenerate real linear
transformations of the phase plane. We denote by $r_q$ the linear presentation
of
 any element $q\in Q$ into the coefficient space $R^{12}$ of system (1).
\smallskip

{\bf Definition 1.} {\rm [15]}  {\sl  A polynomial $K(a,x)$
of the coefficients of system (1) and the  unknown variables  $x^1$ and $x^2$
is called a comitant of system (1) in the group  $Q$, if there exists a
function $\lambda(q)$ such that
$$
K(r_q\cdot a, q\cdot x)\equiv\lambda(q) K(a,x)
$$
for every  $q\in Q,\ a\in R^{12}\; $and $ x = (x^1,x^2)$. }
\smallskip

  A comitant $K$ of system (1) in the group $Q=GL(2,R)$ of linear
homogeneous transformations of the phase plane of system (1)
(which is also called  a group of center-affine transformations) is called
{\sl center affine}.  A comitant $K$ of the
system (1) in the group $Q=Aff(2,R)$
of affine (linear non-homogenous) transformations is called {\sl affine}.
If the comitant $K$ does not depend explicitly on the variables $x^1$ and
$x^2$ then it is called an invariant (center affine or affine, respectively).
\smallskip

{\bf Remark 1.} {\sl
  We say that the comitant of system (1) equals  zero when all its
coefficients vanish. The signs of the comitants which take part in some
sequences of conditions should be calculated  at one and the same point,
where  they do not vanish.}
\smallskip

 We denote by $T(2,R)$ the group of shift transformations and by $r_t$
 the  linear presentation of any element $t\in T$ into the
coefficient space $R^{12}$ of system (1).
\smallskip

{\bf Definition 2.} {\rm [16]} {\sl A comitant $K(a,x)$
of system (1) is called a $\quad T$-comitant if the relation
$$
K(r_t\cdot a,\, x)\, \equiv K(a,\, x)
$$
is valid for every $t\in T$ and $a\in
R^{12}$.}
\smallskip

{\bf Definition 3.} {\rm [17]} {\sl
The polynomial
$$
  (f,\varphi)^{(k)}={(r-k)!(\rho-k)!\over r!\rho!}\sum_{h=0}^k (-1)^h C_k^h
   {\partial^k f\over \partial (x^1)^{k-h}\partial (x^2)^h}
   {\partial^k \varphi\over \partial (x^1)^k\partial (x^2)^{k-h}}
$$
is called a transvectant of index $k$ of two forms $f$ and $\varphi$.
The degree of these forms in the coordinates of the
vector $x=(x^1,x^2)$ are equal to
$r$ and $\rho$, respectively and  $k\le \min(r,\rho)$.}
\smallskip

{\bf Proposition 1.} {\rm [16]} {\sl
The transvectant $(f,\varphi)^{(k)}$ of two $T$-comitants $f$ and $\varphi$
is also  a $T$-comitant.}
\smallskip

 According to [16], by  using  the following $T$-comitants
$$
\eqalign{
   \hat A\, =&\, a^p_k a^q_{\alpha m}a^\alpha_{ln}
           \varepsilon_{pq}\varepsilon^{kl}\varepsilon^{mn}, \cr
   \hat B\, =&\,[2a^n a^h_{u\alpha} - a^n_ua^h_{\alpha}]
    a^l_ra^k_{p\beta}a^m_{qs}a^g_{v\gamma}x^{\alpha}x^{\beta}x^{\gamma}
               \varepsilon_{kl}\varepsilon_{mn}\varepsilon_{gh}
               \varepsilon^{pq}\varepsilon^{rs}\varepsilon^{uv},\cr
  \hat C\, =&\, a^p_{\alpha\beta}x^q x^\alpha x^\beta \varepsilon_{pq},\cr
  \hat D\, =&\,[2a^p a^r_{\alpha \gamma}
                - a^p_\alpha a^r_\gamma ]a^u_{\beta\varkappa}x^q x^s x^v
               \varepsilon_{pq}\varepsilon_{rs}\varepsilon_{uv}
               \varepsilon^{\alpha\beta}\varepsilon^{\gamma\varkappa},\cr
  \hat E\, =&\, a^p_k a^q_{\alpha m}a^r_{ln}x^s x^\alpha\varepsilon_{pq}
                \varepsilon_{rs}\varepsilon^{kl}\varepsilon^{mn},\cr
  \hat F\, =&\,[a^m_sa^n_{\beta}a^k_{pr}-
    2a^k_ra^n_{\beta}a^m_{ps}+
    a^k_pa^m_ra^n_{s\beta}-
    4a^ma^k_{pr}a^n_{s\beta}]a^l_{q\alpha}
    x^{\alpha}x^{\beta}
               \varepsilon_{kl}\varepsilon_{mn}
               \varepsilon^{pq}\varepsilon^{rs},\cr
   \hat G\, =&\, a^\alpha_{\alpha\beta}x^\beta,\cr
   \hat H\, =&\,\frac{1}{2} a^p_{r\alpha}a^q_{s\beta}x^\alpha x^\beta\varepsilon_{pq}
              \varepsilon^{rs},\cr
   \hat K\, =&\, \frac{1}{2} a^p_{mu}a^r_{nv}x^q x^s\varepsilon_{pq}\varepsilon_{rs}
            \varepsilon^{mn}\varepsilon^{uv} \cr
}
$$
(where
$ \varepsilon^{11}= \varepsilon^{22}=\varepsilon_{11}=\varepsilon_{22}=0,\
\varepsilon^{12}=\varepsilon_{12}=-\varepsilon^{21}=-\varepsilon_{21}=1)$,
in virtue of Proposition 1,  the following affine invariants can be
constructed:
$$
\eqalign{
A_{1}  \ =&\ \hat A,\ \
A_{2}  \ = \ (\hat C, \hat D)^{(3)},\ \
A_{3}  \ = \ (((\hat C, \hat G)^{(1)}, \hat G)^{(1)}, \hat G)^{(1)}
,\cr
A_{4}  \ =&\ (\hat H, \hat H)^{(2)},\
A_{5}  \ = \ (\hat H, \hat K)^{(2)},\
A_{6}  \ = \ (\hat E, \hat H)^{(2)}
,\cr
A_{7}  \ =&\ ((\hat C, \hat E)^{(2)}, \hat G)^{(1)},\
A_{8}  \ = \ ((\hat D, \hat H)^{(2)}, \hat G)^{(1)}
,\cr
A_{9}  \ =&\ (((\hat D, \hat G)^{(1)}, \hat G)^{(1)}, \hat G)^{(1)},\
A_{10} \ = \ ((\hat D, \hat K)^{(2)}, \hat G)^{(1)}
,\cr
A_{11} \ =&\ (\hat F, \hat K)^{(2)},\
A_{12} \ = \ (\hat F, \hat H)^{(2)},\
A_{13} \ = \ (((\hat C, \hat H)^{(1)}, \hat H)^{(2)}, \hat G)^{(1)}
,\cr
A_{14} \ =&\ (\hat B, \hat C)^{(3)},\
A_{15} \ = \ (\hat E, \hat F)^{(2)},\
A_{16} \ = \ (((\hat E, \hat G)^{(1)}, \hat C)^{(1)}, \hat K)^{(2)}
,\cr
A_{17} \ =&\ (((\hat D, \hat D)^{(2)}, \hat G)^{(1)}, \hat G)^{(1)},\
A_{18} \ = \ ((\hat D, \hat F)^{(2)}, \hat G)^{(1)}
,\cr
A_{19} \ =&\ ((\hat D, \hat D)^{(2)}, \hat H)^{(2)},\
A_{20} \ = \ ((\hat C, \hat D)^{(2)}, \hat F)^{(2)}
,\cr
A_{21} \ =&\ ((\hat D, \hat D)^{(2)}, \hat K)^{(2)},\
A_{22} \ = \ (((((\hat C, \hat D)^{(1)}, \hat G)^{(1)}, \hat G)^{(1)}, \hat G)^{(1)}  \hat G)^{(1)}
,\cr
A_{23} \ =&\ ((\hat F, \hat H)^{(1)}, \hat K)^{(2)},\
A_{24} \ = \ (((\hat C, \hat D)^{(2)}, \hat K)^{(1)}, \hat H)^{(2)}
,\cr
A_{25} \ =&\ ((\hat D, \hat D)^{(2)}, \hat E)^{(2)},\
A_{26} \ = \ (\hat B, \hat D)^{(3)}.
}
$$
Now, we can introduce the following affine invariants
$$
\eqalign{
C_1 & = 15 A_2^2 - 33 A_{17} -8 A_{18} -63 A_{19}-6 A_{20}-9 A_{21},\cr
C_2 & = -3  A_1A_2 + 2 A_{15},\ C_3= A_2,\ C_4= A_7,\ C_5=- A_2A_3+ 2 A_{22},\cr
C_6 & = 12 A_1(4A_6 -A_7)+3A_2(4 A_4-A_3+12A_5)+6 A_{22}+ 16 A_{23},\cr
C_7 & = -20 A_1A_7 - A_2A_3+ 2 A_{22},\ E_1  =  A_5,\ E_2= A_{25},\cr
C_8 & =  -6 A_1^2 -5 A_8 - A_{10} - A_{11} -3 A_{12},\ C_9=  A_4-A_5,\cr
C_{10} &= A_{26},\
C_{11}=  A_2^2 -10 A_{17} -2 A_{18} -6 A_{19} + 6 A_{21}, \cr
C_{12} & =  - A_1^2(10A_3+9A_5)-3 A_1A_{16}+ A_3(30A_8-7A_{10}-5A_{11})+\cr
  &\quad +A_4(-22A_8+18A_9-11A_{10}+3A_{11})+18A_7(3A_6+5A_7)+\cr
  &\quad +48 A_2A_{13}-2A_6^2+ A_5(46A_8-2A_9+5A_{10}-9A_{11}),
}\eqno(3)
$$
as polynomials in elements $A_1-A_{26}$.

  On the other hand let us consider the following center affine invariants
and comitants, which
are constructed directly through the right-side parts of system (1):
$$
\eqalign{
J_{1}\, = &\,
   \left\vert\,\matrix
               \partial P_1/\partial x^1 &
               \partial P_1/\partial x^2 \cr
               \partial Q_1/\partial x^1 &
               \partial Q_1/\partial x^2
 \endmatrix \right\vert\  = \
 a^\alpha_p a^\beta_q
\varepsilon_{\alpha \beta}
\varepsilon^{pq},
\cr
B_{1}\, = &\,
   \left\vert\,\matrix
               \partial P_1/\partial x^1 &
               \partial P_2/\partial x^2 \cr
               \partial Q_1/\partial x^1 &
               \partial Q_2/\partial x^2
 \endmatrix \right\vert -
   \left\vert\,\matrix
               \partial P_1/\partial x^2 &
               \partial P_2/\partial x^1 \cr
               \partial Q_1/\partial x^2 &
               \partial Q_2/\partial x^1
 \endmatrix \right\vert =
 x^\alpha
 a^\beta_q
a^\gamma_{p \alpha}
\varepsilon_{\beta \gamma}
\varepsilon^{pq},
\cr
B_{2} = &
   \left\vert\,\matrix
                 P_0 &  P_1 \cr
                 Q_0 &  Q_1 \cr
 \endmatrix \right\vert =
 x^\alpha
 a^\beta
 a^\gamma_\alpha
\varepsilon_{\beta \gamma},\
B_{3} =
   {1\over4}\left\vert\,\matrix
               \partial P_2/\partial x^1 &
               \partial P_2/\partial x^2 \cr
               \partial Q_2/\partial x^1 &
               \partial Q_2/\partial x^2 \cr
 \endmatrix \right\vert = \hat H,
\cr
B_{4} = &
   \left\vert\,\matrix
                 P_0 &  P_2 \cr
                 Q_0 &  Q_2 \cr
 \endmatrix \right\vert =
 x^\alpha x^\beta
 a^\gamma
a^\delta_{\alpha \beta}
\varepsilon_{\gamma \delta},\
B_{5} =
   \left\vert\,\matrix
                 P_1 &  P_2 \cr
                 Q_1 &  Q_2 \cr
 \endmatrix \right\vert =
 x^\alpha x^\beta x^\gamma
 a^\delta_\alpha
a^\mu_{\beta \gamma}
\varepsilon_{\delta \mu}.
}
$$

We shall construct the following  transvectants
$$
\eqalign{
   \mu_1&=(B_{3},B_{3})^{(2)},\
   H_1 = (B_{3},B_{1})^{(1)},\
   G_1 = (B_{1},B_{5})^{(1)},\cr
   G_2& = (B_{5},B_{5})^{(2)},\
   G_3 = (B_{3},B_{4})^{(1)},\
   D_1 = \ (((\hat D, \hat D)^ {(2)}, \hat D)^{(1)}, \hat D)^{(3)},
}
$$
   and introduce the following notations
$$
\eqalign{
\mu=& -2\mu_1,\ \
H=2H_1,\ \
2G= 4G_1 -3G_2+ 8G_3,\ \tilde S= B_3,
,\cr
F = & J_1B_5 + 2 B_1B_4+ 4B_2B_3, \
V =  B_4^2 - B_2B_5,\
 \tilde N=\hat K,\
 D= - D_1.
} \eqno(4)
$$
As it was shown in [14], the comitants $\mu,\ H,\ G,\ F$ and $V$ are
responsible for the number and multiplicities of the finite singular
points (FSP) of the quadratic system (1). According to [14] we can construct the
following $T$-comitants:
$$
\eqalign{
 P\,=&\, G^2-6FH+12\mu V,\ R\,=\,4(3H^2-2\mu G),\ S\,=\,R^2-16\mu^2P,\cr
 T\,=&\, 2\mu[2G^3+9\mu(3F^2-8GV)-18FGH+108H^2V]- PR,\; U=F^2-4GV.
}
$$

  We denote by $r_j\,(c_j)$ any real (complex) FSP of the
system (1) of multiplicity $j$ and by $m_f$ [18], we denote the sum of the
multiplicities of all FSP (real and complex) of this system.
\smallskip

{\bf Proposition 2.} {\rm [14]} {\sl The number and multiplicity of the finite
singular points of system (1) are determined in Table 1. The
sets $M_j\in R ^{12}\ (j=1,2,...,19)$ which are defined by the
conditions given in the third and the fourth columns constitute an affine
invariant partition of $R ^{12}$, that is,
$$
\bigcup_{j=1}^{19}M_j=R ^{12},\quad M_i\bigcap_{i\ne j}M_j=\emptyset
$$
and each set $M_j$ is affine invariant.}
\hfill Table 1
\smallskip

\vbox{\tabskip=0pt\offinterlineskip
\halign to 12.785cm {\strut#%
&\vrule#\tabskip=5pt plus 15pt &\hfil#\hfil
&\vrule#&\hfil#\hfil
&\vrule#&\hfil#\hfil
&\vrule#&#\hfil
&\vrule#&#\hfil
&\vrule#\tabskip=0pt\cr
\noalign{\hrule}
&&\multispan3\hfil
Singular points
\hfil&&\multispan3\hfil
Conditions on
\hfil&&&\cr
\noalign{\hrule}
&&\hfil $m_f$ \hfil&&\hfil
multiplicity
\hfil&&\hfil
invariants
\hfil&&\hfil
comitants
\hfil&&\hfil $M_j$ \hfil&\cr
\noalign{\hrule}
&& 4 && $r_1\,r_1\,r_1\,r_1$
&& $\mu\ne0,\,D<0$
&& $R>0,\,S>0$
&& $M_1$ &\cr
&& 4 && $r_1\,r_1\,c_1\,c_1$
&& $\mu\ne0,\,D>0$
&& \hfil -- \hfil
&& $M_2$ &\cr
&& 4 && $c_1\,c_1\,c_1\,c_1$
&& $\mu\ne0,\,D<0$
&& $(R\le0)\vee(S\le0)$
&& $M_3$ &\cr
&& 4 && $r_2\,r_1\,r_1$
&& $\mu\ne0,\,D=0$
&& $T<0$
&& $M_4$ &\cr
&& 4 && $r_2\,c_1\,c_1$
&& $\mu\ne0,\,D=0$
&& $T>0$
&& $M_5$ &\cr
&& 4 && $r_2\,r_2$
&& $\mu\ne0,\,D=0$
&& $T=0,\,PR>0$
&& $M_6$ &\cr
&& 4 && $c_2\,c_2$
&& $\mu\ne0,\,D=0$
&& $T=0,\,PR<0$
&& $M_7$ &\cr
&& 4 && $r_3\,r_1$
&& $\mu\ne0,\,D=0$
&& $T=0,\,P=0,\,R\ne0$
&& $M_8$ &\cr
&& 4 && $r_4$
&& $\mu\ne0,\,D=0$
&& $T=0,\,P=0,\,R=0$
&& $M_9$ &\cr
&& 3 && $r_1\,r_1\,r_1$
&& $\mu=0,\,D<0$
&& $R\ne0$
&& $M_{10}$ &\cr
&& 3 && $r_1\,c_1\,c_1$
&& $\mu=0,\,D>0$
&& $R\ne0$
&& $M_{11}$ &\cr
&& 3 && $r_2\,r_1$
&& $\mu=0,\,D=0$
&& $R\ne0,\,P\ne0$
&& $M_{12}$ &\cr
&& 3 && $r_3$
&& $\mu=0,\,D=0$
&& $R\ne0,\,P=0$
&& $M_{13}$ &\cr
&& 2 && $r_1\,r_1$
&& $\mu=0$
&& $R=0,\,P\ne0,\,U>0$
&& $M_{14}$ &\cr
&& 2 && $c_1\,c_1$
&& $\mu=0$
&& $R=0,\,P\ne0,\,U<0$
&& $M_{15}$ &\cr
&& 2 && $r_2$
&& $\mu=0$
&& $R=0,\,P\ne0,\,U=0$
&& $M_{16}$ &\cr
&& 1 && $r_1$
&& $\mu=0$
&& $R=0,\,P=0,\,U\ne0$
&& $M_{17}$ &\cr
&& 0 && \hfil -- \hfil
&& $\mu=0$
&& $R=P=U=0,\,V\ne0$
&& $M_{18}$ &\cr
&& $\infty$ && \hfil -- \hfil
&& $\mu=0$
&& $R=P=U=0,\,V=0$
&& $M_{19}$ &\cr
\noalign{\hrule}
}}
\smallskip

 Herein we shall use the following assertion from [15] (see Theorem 1.34):
\smallskip

{\bf Proposition 3.} {\sl  For the existence of a center of system (1)
at the origin of coordinates (i.e. $a^1=a^2=0$) it is necessary and sufficient
that the following conditions hold:
$$
   I_1=I_6=0,\quad I_2<0,
$$
and that at least one of the following three  conditions be met:
$$
  1)\ I_3=0;\ 2)\ I_{13}=0;\ 3)\ 5I_3-2I_4\ =\ 13I_3-10I_5=0,
$$
where
$$
\eqalign{
I_{1} =&
 a^\alpha_\alpha,\ \
I_{2} =
 a^\alpha_\beta a^\beta_\alpha,\ \
I_{3} =
 a^\alpha_p
a^\beta_{\alpha q}a^\gamma_{\beta \gamma}
\varepsilon^{pq},\ \
I_{4} =
 a^\alpha_p
a^\beta_{\beta q}a^\gamma_{\alpha \gamma}
\varepsilon^{pq},\cr
I_{5} = &
 a^\alpha_p
a^\beta_{\gamma q}a^\gamma_{\alpha \beta}
\varepsilon^{pq},\
I_{6} =
 a^\alpha_p a^\beta_\gamma
a^\gamma_{\alpha q}a^\delta_{\beta \delta}
\varepsilon^{pq},\
I_{13} =
 a^\alpha_p
a^\beta_{q r}a^\gamma_{\gamma s}a^\delta_{\alpha \beta}a^\mu_{\delta \mu}
\varepsilon^{pq}\varepsilon^{rs}.
}
$$
}
\smallskip
   Herein we have used the notations for invariants from [15].

\head
 Main results
\endhead
\smallskip
First we prove the following lemma:
\smallskip
{\bf Lemma 1.} {\sl  For the existence of a center arbitrarily situated in
the phase
plane of system (1) it is necessary that conditions $C_1=C_3=0$ be
satisfied.}
\smallskip

  Indeed, let system (1) have a singular point of the center type. After
the translation the origin of coordinates to this point we obtain the system
$$
\eqalign{
    {dx\over dt}&= cx + dy + ex^2 + 2hxy + ky^2,\cr
    {dy\over dt}&= ex + fy + lx^2 + 2mxy + ny^2,
}
$$
for which
$ C_1= (c+f)\bar C_1$,
where $\bar C_1$ is a polynomial on the
parameters of this system. Since (according
to Proposition 3) the necessary condition for the existence of a center in (0,0)
is $I_1=c+f=0$ we obtain $C_1=0$. Hereby, for the given system we obtain
$ C_3= I_6$,
and, according to Proposition 3, a center can occur only if $I_6=0,$  i.e. $C_3=0$.

\bigskip
\centerline{\bf \S 1. System with total multiplicity $m_f= 4$}
\bigskip

Herein we shall find out conditions for the existence of a center by using
invariants (3) and Table 1 in the case where the
multiplicity $m_f$ equals four.
From Table 1 and [19] it follows that the
system (1) with $m_f=4$ may have a center
only if it belongs to the set $M_1\cup M_2\cup M_4\cup M_8$. This implies 4 different
cases which will be examined in the sequel.
\smallskip

{\bf Theorem 1.} {\sl System (1) with conditions $\mu\ne0,\ D<0,\ R>0,
\ S>0$  (there are 4 simple singular points) has one center if and  only if
one of the following two  sequences of conditions holds:
$$
\eqalign{
   (i)\ &\ C_2C_4<0 ,\ C_1=C_3=C_5=0;\cr
   (ii)\,&\ C_4=0,\  C_1=C_3=0,\  \mu<0;\cr
}
$$
and it has two centers if and only if the following conditions hold}
$$
\eqalign{
   (iii)&\ C_4=0,\  C_1=C_3=0,\   C_9\ge0,\ \mu>0.\cr
}
$$
\smallskip

{\bf Proof.} Let us assume that system (1) has four  real distinct
singular points.  By applying an  affine transformation we can assume that
three  singular
such points are the points $M_0(0,0),\ M_1(1,0)$ and $M_2(0,1)$, respectively.
Hence, system (1) becomes
$$
\eqalign{
    {dx\over dt}&= cx + dy - cx^2 + 2hxy - dy^2,\cr
    {dy\over dt}&= ex + fy - ex^2 + 2mxy - fy^2,
}\eqno(5)
$$
from which, by using the notations
$\bar C= cm-eh,\ \bar D=de-cf$ and $\bar F=fh-dm$,
we obtain
$$
  \mu= \bar D^2 - 4 \bar C\bar F,\
  D=-\bar D^2(\bar D-2\bar C)^2(\bar D-2\bar F)^2.
$$
   It is easy to observe, that system (5),  besides the critical points
$ M_0(0,0)$, $ M_1(1,0)$ and $ M_2(0,1)$,
has the critical point $ M_3(x_0,y_0)$
with coordinates
$$
   x_0= \bar D(\bar D-2\bar F)/\mu, \quad
   y_0= \bar D(\bar D-2\bar C)/\mu.
$$
  Following  [19], we shall calculate the coefficients of the
equation
$$
   \lambda^2 - \sigma \lambda + \Delta =0,
$$
which determine the eigenvalues corresponding to a critical point.
Thus, for the critical points $M_0,\ M_1,\ M_2$ and $M_3 $ we obtain
$$
\eqalign{
  \sigma^{(0)}&=c+f,\ \Delta^{(0)}=-\bar D,\cr
  \sigma^{(1)}&=-c+f+2m,\ \Delta^{(1)}=\bar D-2\bar C,\cr
  \sigma^{(2)}&= c-f+2h,\ \Delta^{(2)}=\bar D-2\bar F,\cr
  \sigma^{(3)}&= c+f + 2(m-c)x_0 + 2(h-f)y_0 ,\cr
   \Delta^{(3)}&=-\bar D(\bar D-2\bar C)(\bar D-2\bar F)/\mu,\
} \eqno(6)
$$
respectively.

 For system (5) we obtain
$$
     C_1=3\mu\sigma^{(0)}\sigma^{(1)}\sigma^{(2)}\sigma^{(4)}
$$
and, according to Lemma 1,  if $C_1\ne0$ system (5) has not a center.

    If $C_1=0$,  in virtue of $\mu\ne0$ it follows  that
$\sigma^{(0)}\sigma^{(1)}\sigma^{(2)}\sigma^{(3)}=0$.
Without loss of generality we can suppose that $\sigma^{(0)}=0$, otherwise
we can use the linear transformation, which replaces the points $M_0$ and
$M_i$  in  case $\sigma^{(i)}=0 (i=1,2,3)$. In all three cases, after the
corresponding change of parameters we obtain  the same system, but with
$\sigma^{(0)}=0$.

Thus, from $\sigma^{(0)}=0$ and (6) we get $f=-c$ and for system (5)
we have
$$
   C_4= \mu\sigma^{(1)}\sigma^{(2)}\sigma^{(3)}/\bar D. \eqno(7)
$$

  For the
  system (5) with the
  condition $f=-c$, the following affine invariants can be
calculated:
$$
\eqalign{
  & C_4={1\over 3}I_{13}^{(0)},\ C_2=3I_2^{(0)}C_4,\
         C_3=-{2\over 3}I_{6}^{(0)}, \cr
  & C_5=6I_3^{(0)}C_4,\ C_6=6(5I_3^{(0)}-2I_4^{(0)})C_4,\
         C_7=2(13I_3^{(0)}-10I_5^{(0)})C_4\cr
}\eqno(8)
$$
where $I_j^{(0)}\ (j=2,3,4,5,6,13)$ are the values of the center affine
invariants from Proposition 3, calculated for system (2) with the
singular
point $M_0(0,0)$.
\smallskip
  I. If $C_4\ne0$, in accordance with (7) we obtain
$\sigma^{(1)}\sigma^{(2)}\sigma^{(3)}\ne0$ and neither of the singular points
$M_i (i=1,2,3)$ can be a center. In this case, taking into account
Proposition 3, the
relation $I_1=c+f=0$ and (8), we conclude, that the singular point
$M_0(0,0)$ will be a center if and only if the following conditions hold:
$$
     C_1=C_3=C_5(C_6^2+C_7^2)=0,\quad C_2C_4<0.
$$
However, we shall prove that in the case under consideration conditions
$C_5\ne0$ and $C_6=C_7=0$ can not be satisfied.  Indeed, if we suppose
the contrary, then from (8) it follows:
$ I_{13}\ne0,\ I_2<0, \ I_1=I_6=0$, $I_3\ne0$ and $5I_3-2I_4=13I_3-10I_5=0$.

Hereby, as it was shown in [20, p. 131], by applying a center affine
transformation  system (1) can be brought to the canonical system
$$
   {dx\over dt} =y+ qx^2+ xy,\quad {dy\over dt}=-x - x^2+3qxy + 2y^2,
$$
for which $D=8q^2(q^2+1)>0$ in virtue of $I_{13}=125q(q^2+1)/8\ne0$.
As we can observe this contradicts condition $D<0$ of Theorem 1.
\smallskip

  II. Let condition $C_4=0$ be satisfied. From (7) we obtain that
$\sigma^{(1)}\sigma^{(2)}\sigma^{(3)}=0$ and without loss of generality, we
can assume that $\sigma^{(1)}=0$, otherwise we can use a linear
transformation. Indeed, if $\sigma^{(2)}=0$ (resp. $\sigma^{(3)}=0$)
the transformation
$x_1=y,\ y_1=x \ (x_1=x/x_0,\ y_1=y-xy_0/x_0)$  replaces the points
$M_1$ and $M_2$ (resp. $M_1$ and $M_3$ ) and keeps the other points.

   Thus, the conditions $C_1=C_4=0$ imply that $\sigma^{(0)}= \sigma^{(1)}=0$ and from
(6) we receive
$$
   f=-c,\ m=c,\ \sigma^{(2)}=2(c+h),\ \sigma^{(3)}= 2(c+h)y_0  \eqno(9)
$$
and, hence, system (5) becomes as
$$
\eqalign{
    {dx\over dt}&= cx + dy - cx^2 + 2hxy - dy^2,\cr
    {dy\over dt}&= ex - cy - ex^2 + 2cxy + cy^2.
}\eqno(10)
$$
For system (10) the following comitants can be calculated:
$$
\eqalign{
  &  I_1^{(0)}=I_{13}^{(0)}=0,\ I_2^{(0)}=2(c^2+de),\cr
  & I_6^{(0)}= (c+h)(2ceh-2c^3-c^2e-de^2).
}
$$
On the other hand by translating the origin of coordinates at the singular
point $M_1(1,0)$ we obtain the system
$$
\eqalign{
    {dx\over dt}&= cx + dy - cx^2 + 2hxy - dy^2,\cr
    {dy\over dt}&= ex - cy - ex^2 + 2cxy + cy^2,
}\eqno(11)
$$
for which
$$
\eqalign{
  &  I_1^{(1)}=I_{13}^{(1)}=0,\ I_2^{(1)}=2(c^2-de-2eh),\cr
  & I_6^{(1)}= (c+h)(2ceh-2c^3-c^2e-de^2).
}
$$
 For system (10), as well as for system (11) we obtain
$$
\eqalign{
 & C_3=-{2\over 3}I_{6}^{(0)}=-{2\over 3}I_{6}^{(1)},\cr
 & C_8=  -{4\over 3}(c+h)^2(c^2+de)(c^2-de-2eh)=
 -{1\over 3}(c+h)^2I_{2}^{(0)}I_{2}^{(1)},\cr
 & C_9=-(c+h)^2(c^2-eh) =-{1\over 2}(c+h)^2(I_{2}^{(0)}+I_{2}^{(1)}), \cr
 & D=-{1\over 18}(c^2+de)^2(c^2-de-2eh)^2(c^2+de+2cd+2ch)^2\ne0.\cr
}\eqno(12)
$$

  1) If $C_8\ne0$, from (12) we obtain $(c+h)\ne0$ and from (9) we
get
$\sigma^{(2)}\sigma^{(3)}\ne0$. Therefore, from (12) and
Proposition 3, we conclude, that for $ C_8\ne0$ system (1) has
one center if and only if $C_3=0$ and  $ C_8>0$
and has two center if and only if $C_3=0,\ C_8<0$ and $C_9>0$.

 It is not too difficult to show that conditions $C_3=0$ and  $ C_8\ne0$ imply
 $\mu C_8<0$. Indeed, by virtue of (12) from $C_3=0$ and
 $ C_8\ne0$  results $e\ne0$. Therefore we can consider $e=1$ and, hence,
by (12) condition $C_3=0$ yields $d=2ch-2c^3-c^2$. Hereby for system (11)
we have
$$
   \mu=-4c(c+1)(c^2-h)^2,\ C_8=\frac{16}{3}c(c+1)(c+h)^2(c^2-h)^2,
$$
and this implies $\mu C_8<0$.
\smallskip

  2) Let us assume that the condition $C_8=0$ holds.  Since condition
$D\ne0$ is satisfied, from (12) we get $(c+h)=0$, and from
(9) and (6) we obtain
$$
   \sigma^{(0)}=\sigma^{(1)}=\sigma^{(2)}= \sigma^{(3)}=0,\
   Sgn(\Delta^{(0)}\Delta^{(1)}\Delta^{(2)}\Delta^{(3)})=Sgn\mu.
\eqno(13)
$$
On the other hand, following
Proposition 3 we  calculate the values of the invariants $I_1,\ I_2,\ I_6$
and $I_{13}$ for the
system (5) and for other three systems, that are obtained
from (5) by
placing the critical points $M_1,\ M_2$ and $M_3$ at the origin, respectively.
Thus, we get the following expressions:
$$
\eqalign{
   & I_{1}^{(i)}=I_{6}^{(i)}=I_{13}^{(i)} =0\ (i=0,1,2,3),\cr
   & I_{2}^{(0)}=-2\Delta^{(0)},\ I_{2}^{(1)}=-2\Delta^{(1)},\cr
   & I_{2}^{(2)}=-2\Delta^{(2)},\  I_{2}^{(3)}=-2\Delta^{(3)}.
}\eqno(14)
$$
 According to Proposition 3, the number of the negative quantities
among the $ I_{2}^{(i)}\ (i=0,1,2,3)$ coincide with the number of  centers
of system (5).

  As it is well known for quadratic system  at least one of the
quantities $\Delta^{(i)}\ (i=0,1,2,3)$ is negative and at least one of them
is positive. On the other hand, at most two singular points can be of the
center or of the focus type. Therefore, from  Proposition 3,  (13) and (14)
we conclude, that system (5) has two centers
if $\mu>0$ and one center if $\mu<0$.
It remain to notice, that
condition $ C_8=0$ (i.e. $c+h=0$) is equivalent to $ C_9=0$. Indeed, if
we suppose that $c+h\ne0$, then by (12) condition $ C_9=0$ yields
$c^2=eh$ and from $C_3=0$ it results $e(c^2+de)=0$, contrary to the
condition  $D\ne0$.

Theorem 1 is proved.
\bigskip

 {\bf Theorem 2.} {\sl System (1) with conditions $\mu\ne0,\ D>0 $
 (there are 2 simple real and two imaginary singular points) has one center
if and  only if one of the following two sequences of conditions holds:
$$
\eqalign{
   (i)\ &\ C_2C_4<0 ,\ C_1=C_3=C_5(C_6^2+C_7^2)=0;\cr
   (ii)\,&\ C_4=0,\ C_{12}\le0,\ C_1=C_3=0,\  \mu>0;\cr
}
$$
and it has two centers if and only if the following  sequence
of conditions holds}
$$
   (iii)\  C_4=0,\ C_{12}<0,\ C_1=C_3=0,\  \mu<0,\ C_9>0.
$$
\bigskip
{\bf Proof.} Let us assume that system (1) has two  real and two imaginary
singular points.  By applying the affine transformation we can move two
real  singular  points to the points $M_0(0,0)$ and $M_1(1,0)$, respectively.
Hence, system (1) becomes
$$
\eqalign{
    {dx\over dt}&= cx + dy - cx^2 + 2hxy + ky^2,\cr
    {dy\over dt}&= ex + fy - ex^2 + 2mxy + ny^2,
}\eqno(14)
$$
for which, by using the notations
$$
\eqalign{
  \bar B= cn-ek,\ \bar C= cm-eh\ \bar D=de-cf,\cr
  \bar E= dn-fk,\ \bar F=fh-dm,\ \bar H=hn-km
}
$$
 we obtain
$$
\eqalign{
 & \mu= \bar B^2 + 4 \bar C\bar H\ne0,\ D=-{2\over 27}\bar D^2(\bar D-2\bar C)^2 Z>0,\cr
 &  Z=(\bar B-2\bar F)^2 + 4\bar E(\bar D-2\bar C).
}\eqno(15)
$$
Since  $D>0$, it follows that $Z<0$ and for the
real singular points $M_0(0,0)$
and $M_1(1,0))$ we obtain
$$
  \sigma^{(0)}=c+f,\ \sigma^{(1)}=-c+f+2m. \eqno(16)
$$
For system (14)  we can  calculate that
$$
\eqalign{
\mu C_1 & = 3 \sigma^{(0)}\sigma^{(1)}(R^2-ZI^2),\
 \mu  C_4  =\bar P R+\bar Q I,\cr
 C_{12}& = {1\over 4}[\mu^2(\sigma^{(0)}-\sigma^{(1)})^2-
   4mR(\sigma^{(0)}+\sigma^{(1)})+4ZI^2],\
}\eqno(17)
$$
where
$$
\eqalign{
  R&= (c-m)(2\bar B\bar F +4\bar C\bar E-\bar B^2)
       -2(h+n)(\bar B\bar C -\bar B\bar D + 2\bar C\bar F)+\cr
     &\ \ +(c+f)\mu,\ I= (m-c)\bar B -2(h+n)\bar C,\cr
\bar P&=\bar B(-cn-hm+mn)- \bar C[(h-n)^2+k(c-2m)]+\bar Hm^2,\cr
\bar Q&=\bar B^2(2fh-cd-ch+2hm) +\bar B \bar D(-ck-h^2+2hn-n^2)+\cr
    &\ \   +\bar B \bar C(ck-4dh-2fk-3h^2-2hn-2km+ n^2)
           - \bar B \bar Fm(2h+n)-\cr
    &\ \   -\bar B \bar Hm(f+m) +2\bar C^2k(d+2h+2n)-2\bar C \bar Dk(h+2n)+\cr
    &\ \    +2\bar C \bar F(h+n)^2
           +4\bar C \bar Hm(d+h+n) -4\bar D \bar Hmn +2\bar F \bar Hm^2.\cr
}\eqno(18)
$$

 According to Lemma 2,  if $C_1\ne0$  system (14) has not a center.

   Let us assume that condition $C_1=0$ holds.
\smallskip

{\bf Case I.} If $C_4\ne0$ from (18) we get $R^2+I^2\ne 0$ and
by $Z<0$ and (17) the condition $C_1=0$ implies that
$\sigma^{(0)}\sigma^{(1)}=0$.
Without loss of generality we can assume that $\sigma^{(0)}=0$, otherwise
we can use the linear transformation, which  replaces the  points $M_0$ and
$M_1$.

Thus,  $\sigma^{(0)}=0$, and from (16), it follows that
 $f=-c$ and for system (14)
we obtain
$$
   12\mu\bar DC_4= \sigma^{(1)}(R^2-ZI^2). \eqno(19)
$$
On the other hand, for system (14) with the
condition $f=-c$, the following affine
invariants can be calculated:
$$
\eqalign{
   & C_4={1\over 3}I_{13}^{(0)},\  C_2=3I_2^{(0)}C_4,\
          C_3=-{2\over 3}I_{6}^{(0)}, \cr
   & C_5=6I_3^{(0)}C_4,\  C_6=6(5I_3^{(0)}-2I_4^{(0)})C_4,\
           C_7=2(13I_3^{(0)}-10I_5^{(0)})C_4,\cr
}\eqno(20)
$$
where $I_j^{(0)}\ (j=2,3,4,5,6,13)$ are the values of the center affine
invariants from Proposition 3, calculated for system (14) with singular
point $M_0(0,0)$.
\smallskip
   Since  $\mu D\ne0$, by (19) condition $C_4\ne0$
implies that $\sigma^{(1)}\ne0$, i.e. singular point $M_1 $ is not a  center.

Thus, in the case of $C_4\ne0$, from
Proposition 3 and relationship $I_1=c+f=0$ and (20) we conclude, that the
 system (14)
has one center if and only if the following conditions hold:
$$
     C_1=C_3=C_5(C_6^2+C_7^2)=0,\quad C_2C_4<0.
$$

\medskip
{\bf Case II.} Let condition $C_4=0$ be satisfied.
\smallskip

 {\bf A)} If $C_8\ne0$ we shall examine two cases: $C_{12}<0$ and $C_{12}\ge0$.
\smallskip
   1) Let us consider firstly that condition $C_{12}<0$ holds. Then
 $R^2+I^2\ne0$, otherwise from  (17) we get
$C_{12}=\mu^2(\sigma^{(0)}-\sigma^{(1)})^2\ge0$.
Therefore, $R^2-ZI^2\ne0$ and
conditions $C_1=C_4=0$, (17) and (19) imply that
$\sigma^{(0)}=\sigma^{(1)}=0$. From (16) it follows that $f=-c,\ m=c$ and the
system (14) becomes as
$$
\eqalign{
    {dx\over dt}&= cx + dy - cx^2 + 2hxy + ky^2,\cr
    {dy\over dt}&= ex - cy - ex^2 + 2cxy + ny^2.
}\eqno(21)
$$
For system (21) the following invariants can be calculated:
$$
\eqalign{
  &  I_1^{(0)}=I_{13}^{(0)}=0,\ I_2^{(0)}=2(c^2+de),\cr
  & I_6^{(0)}= (h+n)(2ceh-2c^3-cen+e^2k).
}
$$
On the other hand, by translating the origin of coordinates to the singular
point $M_1(1,0)$ of system (21) we obtain the system
$$
\eqalign{
    {dx\over dt}&= -cx + (d+2h)y - cx^2 + 2hxy + ky^2,\cr
    {dy\over dt}&= -ex + cy - ex^2 + 2cxy + ny^2,
}\eqno(22)
$$
for which
$$
\eqalign{
  &  I_1^{(1)}=I_{13}^{(1)}=0,\ I_2^{(1)}=2(c^2-de-2eh),\cr
  & I_6^{(1)}= (h+n)(2ceh-2c^3-cen+e^2k).
}
$$
 For systems (21) and (22) we obtain
$$
\eqalign{
 & C_3=-{2\over 3}I_{6}^{(0)}=-{2\over 3}I_{6}^{(1)},  \cr
 & C_8=-{4\over 3}(h+n)^2(c^2+de)(c^2-de-2eh)=
       -{1\over 3}(h+n)^2I_{2}^{(0)}I_{2}^{(1)},\cr
 & C_9=(h+n)^2(eh-c^2)=-{1\over 2}(h+n)^2(I_{2}^{(0)}+I_{2}^{(1)}), \cr
 &  D=-{2\over 27}(c^2+de)^2(c^2-de-2eh)^2 Z \ne0,\cr
 &  C_{12}= I^2Z=4(h+n)^2(c^2-eh)^2Z.
}\eqno(23)
$$
\smallskip
     From  $C_{12}\ne0$ and (23) we get $(h+n)\ne0$. Therefore,
from (23) and Proposition 3, we conclude that for $C_1=C_4=0$ and
$C_{12}<0$  system (1) has one center if and only if $C_3=0$ and $C_8>0$
and has two centers if and only if $C_3=0,\ C_8<0$ and $C_9>0$. It
is not difficult to convince, that in the case under consideration
condition  $C_8<0$ ($C_8>0$) is equivalent to $\mu<0$ ($\mu>0$).
Indeed, condition $C_4=0$ by (20) yields $I_{13}=0$ and it was shown in
[20] that by applying a linear transformation system
(1) with can be brought either to the canonical system [20, p. 103]
$$
   {dx\over dt} =y +2(1-c)xy,\quad {dy\over dt}=-x + dx^2+ cy^2.\eqno(24)
$$
or to the canonical system [20, p. 80]
$$
   {dx\over dt} =-y-cx^2-ay^2,\quad {dy\over dt}=x+bx^2+2cxy,\eqno(25)
$$

For system (24) one can be calculated
$$
\eqalign{
 &\mu=4cd(c-1)^2,\  D=\frac{1}{3}c(d+2-2c)^3,\ C_{12}=-2<0,\cr
 &C_8=4d(d+2-2c),\ C_1=C_3=0.
}
$$
Hereby by virtue of $D>0$ it results $\mu C_8>0$ and since $C_{12}<0$
this has proved our affirmation.

As regards system (25) we obtain
$$
  \mu=a^2b^2+4ac^3,\ D=\frac{1}{6}(a-2c)^2(4ac-b^2-8c^2),\
  C_1=C_3=C_8=C_{12}=0,
$$
i.e. condition $C_8\ne0$ is not satisfied.

\smallskip

  2) Let condition  $C_{12}\ge0$ holds. We shall demonstrate that by virtue
$C_8\ne0$ there is no  center on the phase plane of  system (14).
  Indeed, according to (17), the condition $C_1=0$ implies that
$\sigma^{(0)}\sigma^{(1)}(R^2-ZI^2)=0$.

  If $R^2+I^2\ne0$ then from $C_1=C_4=0,\ Z<0$ and
(19) we get  $\sigma^{(0)}=\sigma^{(1)}=0$, i.e. $f=-c,\ m=c$. Therefore,
for system (14) we obtain
$$
\eqalign{
 & C_{12}= ZI^2,\ C_3={2\over 3}(h+n)(2c^3-2ceh+cen-e^2k),
 \ R=3(d+h)C_3,\
}
$$
and conditions $C_{12}\ge0,\ Z<0$ and $\ R^2+I^2\ne0$ imply $I=0,\ R\ne0$.
Hereby, we get that
$C_3\ne0$ and by Lemma 1 system (14) has no
singular point of the  center type. We note that, according to (23) and
$D\ne0$, in this case  $C_8\ne0$.

   Let conditions  $R=I=0$ be satisfied. Without loss of generability
we can assume that
relation $ce=0$ holds. Indeed, if $e\ne0$, by
applying, the transformation
$x_1=x-cy/e,\ y_1=y$, which keeps the singular points $M_0(0,0)$ and
$M_1(1,0)$ for system (14) the relation $c=0$ will be satisfied. Thus, we
shall consider two cases: $e\ne0, c=0$ and $e=0$.
\smallskip

    a) If $e\ne0, c=0$ for system (14) we have $I= e[2h^2+2hn-km]=0$.
\smallskip

 $\alpha$) If $m\ne 0$ by using the change of the time the condition $m=1$
will be satisfied. Therefore, from $I=0$, we obtain $k=2h^2+2hn$ and for system (14)
in this case one can derive
$$
   R= 4eh(fh-d)[e(h+n)^2+2h+n],\ \mu= 4eh^2[e(h+n)^2+2h+n].
$$
Since $R=0,\ \mu\ne 0$ we get $d=fh$. Thus, for system (14)
we obtain
$$
\eqalign{
&   C_3= -{4\over 3}he[e(h+n)^2-f(f+2)(2h+n)],\ Z=-3hC_3,\cr
 &  C_8=  {1\over 3}f(f+2)\mu ,\
   D=-{2\over 27}e^4h^4f^2(f+2)^2 Z,
}\eqno(26)
$$
and by $Z<0$ it follows $C_3\ne 0$ and hence, according to
Lemma 1 system (14) has not any center. Notice, that in this case from (26)
and $D\ne0$ we get $C_8\ne0$.

\smallskip
 $\beta$) If $m=0$ since $c=0$ for system (14), we have
$$
I=2eh(h+n)=0,\ C_8={4\over 3}[3fR+de^2(d+2h)(h+n)^2],\eqno(27)
$$
and from $I=R=0,\ C_8\ne0$ and (27) it follows  $h=0,\ n\ne0$. Therefore,
for system (14) with $c=m=h=0, ne\ne0$ we obtain:
$$
   R=e^2k(fk-2dn)=0,\ \mu=e^2k^2\ne0.
$$

Thus, the condition $k\ne0$ is satisfied and we can consider $k=1$
by a change of scale
 if necessary. Condition $R=0$ implies $f=2dn$
and for system (14) the following polynomials can be calculated:
$$
   C_3= {2\over 3}en(4d^2n-e),\ Z=-e(4d^2n-e).
$$
Hereby, from $Z<0$ and $ne\ne0$, it follows that $C_3\ne 0$ and
from Lemma 1 we can
conclude that there is not any center on the phase plane of system (14) in
this case under consideration.

\smallskip
   b) If $e=0$ for system (14) we have
$$
\mu = c(cn^2+4hmn-4km^2),   \ I= -c(cn+2hm+mn),
$$
hence, the conditions $I=0,\ \mu\ne0$ imply  $m\ne0,\ cn+2hm+mn=0$.
 As it was mentioned above, we can assume
$m=1$  hence, $h=-{1\over 2}n(c+1)$. Therefore, for system (14) we obtain
$$
   \mu=-c(cn^2+4k+2n^2),\ R=c(f+1)\mu,
$$
and from conditions $\mu\ne0,\ R=0$ we get $f=-1$. Hereby for system (14)
the following relation can be established
$$
 \sigma^{(0)}=-\sigma^{(1)}=c-1, \ C_8=-{1\over 3}(c-1)^2\mu\ne0.
\eqno(28)
$$
From $C_8\ne0$ it follows that
$\sigma^{(0)}\sigma^{(1)}\ne0$ and, hence,
there is no singular point of the center type for system (14).
\bigskip

  {\bf B)}  Let us now assume that the condition $C_8=0$ is satisfied.
  As it was indicated
above from $C_1=C_4=0$ and $R^2+I^2\ne 0$ it follows that $C_8\ne0$. Hence,
for system (14) we obtain $R=I=0$ and we again shall examine two cases:
 $e\ne0, c=0$ and $e=0$.
\smallskip
  a) Let the conditions $e\ne0,\ c=0$ hold. If $m\ne0$ it was shown above
that from  (26) and $D\mu\ne0$ it follows that $C_8\ne0$. Hence, $m=0$, and
taking into account (27) and the relations  $R=I=0,\ C_8=0$ for system (14)
we obtain:
$$
I=2eh(h+n)=0,\ C_8={4\over 3}de^2(d+2h)(h+n)^2,\ \mu= e(ek^2-4h^2n).
$$
From $I=0,\ C_8=0$ we get $h=-n\ne0$. Indeed, if we suppose
that condition $h+n\ne0$ is satisfied, then from $I=C_8=0$ it follows  that
$h=0,\ n\ne0,\ de=0$ and we obtain a contradiction with condition
$D=-{2\over 27}d^4e_4 Z\ne0$. Thus, $h+n=0$  and for the system
(14) with $c=m=0,\ n=-h\ne0$ one can calculate:
$$
   \mu=e(ek^2+4h^3)\ne0, \ C_8= 2f^2\mu.
$$
Hereby, conditions $C_8=0,\ \mu\ne0$ imply $f=0$ and the
system (14) becomes as
system
$$
    {dx\over dt} =  dy  + 2hxy + ky^2,\quad
    {dy\over dt} = ex  - ex^2  - hy^2,
\eqno(29)
$$
for which
$$
   C_1=C_3=C_{12}=0,\  Z= e(ek^2-4d^2h-8dh^2)<0,\ \mu=e(ek^2+4h^3)\ne0.
   \eqno(30)
$$
For system (29) the following invariants (from Proposition 3)
can be calculated.
$$
    I_1^{(0)}=I_{6}^{(0)}=I_{13}^{(0)}=0,\ I_2^{(0)}=2de.
$$
On the other hand by translating the origin of coordinates at the singular
point $M_1(1,0)$ of system (29) we obtain the system
$$
\eqalign{
    {dx\over dt}&=  (d+2h)y + 2hxy + ky^2,\quad
    {dy\over dt} = -ex - ex^2 -hy^2,
}\eqno(31)
$$
for which
$$
\eqalign{
  &  I_1^{(1)}=I_{6}^{(1)}=I_{13}^{(1)}=0,\ I_2^{(1)}=-2(de+2eh).\cr
}
$$
It is easy to observe, that
$$
  I_2^{(0)}I_2^{(1)}=-4de^2(d+2h),\ I_2^{(0)}+I_2^{(1)}=-4eh, \
  Sgn\mu=Sgn(eh). \eqno(32)
$$
Indeed, from (30) it follows  that
$$
   \mu-Z=4eh(d+h)^2,
$$
and since $Z<0$ the condition $\mu>0$ implies $eh>0$. On the other hand
from $\mu=e^2k^2+4eh^3<0$ it results $eh<0$. Therefore, taking into account
condition $Z=e^2k^2-4deh(d+2h)<0$ we obtain that $deh(d+2h)>0$ and, hence,
from (32) it follows that
$$
   Sgn[\mu I_2^{(0)}I_2^{(1)}]=-Sgn[de^3h(d+2h)]=-1,\
   Sgn\mu=-Sgn(I_2^{(0)}+I_2^{(1)}).
$$
Hereby, we have obtained, that for $\mu>0$ it results $I_2^{(0)}I_2^{(1)}<0$
and, according to Lemma 1, either singular point $M_0$ or $M_1$ is of the
center type. If $\mu<0$ then $I_2^{(0)}I_2^{(1)}>0$ and
$I_2^{(0)}+I_2^{(1)}>0$. Thus, both quantities $I_2^{(0)}$ and $I_2^{(1)}$
are positive and by Lemma 1 there is no  center on the phase
plane of system (29).
\smallskip
  b) Let us assume
  that condition $e=0$ hold. As it was mentioned above, by virtue
of conditions $\mu\ne0,\ I=R=0$ we get that  $m=1,\ h=-{1\over 2}n(c+1)$
and $f=-1$. Hereby,
we obtain (28) for system (14) and condition $C_8=0$
in this case implies that $c=1$. Therefore, we obtain the system
$$
\eqalign{
    {dx\over dt}&=  x + dy -  x^2 + 4xy + ky^2,\cr
    {dy\over dt}&= - y + 2xy - 2y^2,
}\eqno(33)
$$
for which
$$
  \mu-Z=-4(d-2n)^2\le0,\ C_{12}=0.
$$
Hence, $\mu<0$ because  $Z<0$. Thus, for system (33) as well as
for the translated system (with singular point $M_1$ at the origin of
coordinates) we obtain that $I_2^{(0)}=I_2^{(1)}=2>0$. By
Proposition 3 system (33) has no a center.

 Thus we have found out, that in the case $C_1=C_3=C_4=C_8=0$ system (14) have
a center if and only if the
condition $\mu>0$ hold. Moreover, the center is unique. It remains to note,
that for $C_1=C_3=C_4=0$ and $\mu D\ne0$ conditions $C_8=0$ and $C_{12}=0$
are equivalent, since for system (24) we have $C_8C_{12}\ne0$ and for
system (25) the conditions $C_8=C_{12}=0$ hold.

   As all cases were examined, Theorem 2 is proved.

\bigskip

{\bf Theorem 3.} {\sl For the existence of a center of system (1) with
$\mu\ne0$, $D=0$, $ T<0 $ (there are two simple and one double singular points)
it is necessary and sufficient that one of the following two sequences of
conditions holds:
$$
\eqalign{
   (i)\, &\  C_2C_4<0 ,\ C_1=C_3=C_5=0;\cr
   (ii)&\ C_4=0,\ \mu>0, \ C_1=C_3=C_8=0.
}
$$
Moreover, the center is unique.}
\smallskip

{\bf Proof.}  According to Proposition 2 if conditions $ \mu\ne0,\ D=0,\ T<0 $
are satisfied the
system (1) has two simple and one double singular points situated
on its phase plane. We shall find out the canonical form of a system (1) with
such points.
 By applying the affine transformation we can move these three  singular
points to the points $M_0(0,0),\ M_1(1,0)$ and $M_2(0,1)$, respectively.
Hence, system (1) can be brought to the system
$$
\eqalign{
    {dx\over dt}&= cx + dy - cx^2 + 2hxy - dy^2,\cr
    {dy\over dt}&= ex + fy - ex^2 + 2mxy - fy^2,
}\eqno(34)
$$
for which, by using the notations $\hat C= cm-eh,\ \hat D=de-cf$ and $
\hat F=fh-dm$, we have
$$
  \mu= \hat D^2 - 4 \hat C\hat F,\ D=-\hat D^2(\hat D-2\hat C)^2
(\hat D-2\hat F)^2.
$$
   It is easy to observe, that system (34) has,
   besides the critical points
$ M_0(0,0)$, $M_1(1,0)$ and $ M_2(0,1)$  the critical point $ M_3(x_0,y_0)$,
where
$$
   x_0= {1\over \mu}\hat D(\hat D-2\hat F), \quad
   y_0= {1\over \mu}\hat D(\hat D-2\hat C).
$$
Thus, we can conclude that in virtue of $D=0$ the point $M_3$ will coincide
with $M_0$ for $\hat D=0$, with point $M_1$ for $\hat D=2\hat C$ and with
$M_2$ for $\hat D=2\hat F$.
Without loss of generality we can assume that the condition $D=0$ implies
that
$\hat D=de-cf=0$ and that the
singular point $M_0$ becomes degenerated (we can removed
the respective points if it is necessary). Therefore, without loss of
generality one can sets
$d=cu,\ f=eu$ and  system (34) becomes
$$
\eqalign{
    {dx\over dt}&= cx + cuy - cx^2 + 2hxy - cuy^2,\cr
    {dy\over dt}&= ex + euy - ex^2 + 2mxy - euy^2.
}\eqno(35)
$$
  By [19], for the critical points $ M_1$ and $M_2 $ of system (35)
we must have
$$
\eqalign{
  \sigma^{(1)}&=-c+eu +2m,\quad
  \sigma^{(2)}=c-eu+2h,
} \eqno(36)
$$
respectively.

 For system (35) we obtain
$$
     \mu = 4u(cm-eh)^2\ne0,\
     C_1=3\mu(c+eu)^2\sigma^{(1)}\sigma^{(2)}.\eqno(37)
$$
If $C_1^2+C_3^2\ne0$   Lemma 1 implies that system (35) has no center.

    Let us assume $C_1=C_3=0$ and examine two cases: $C_4\ne0$ and $C_4=0$.
\bigskip

     1) If $C_4\ne0$ then the
     condition $c+eu\ne0$ holds. Indeed, if $c=-eu$ for
system (35) we obtain that
$$
\eqalign{
      \mu& = 4e^2u(h+mu)^2\ne0, \ C_1=0,\
      C_3={4\over 3}e^2u(h+mu)(h-eu^2 -eu- mu),\cr
    C_4&={2\over 3}e(h+mu)(h-eu^2 -eu- mu)(ue+m)(h-ue),
} \eqno(38)
$$
and from $\mu\ne0$ and $C_3=0$ we obtain that $C_4=0$.

  Thus, $c+eu\ne0$ and from (37) the condition $C_1=0$ implies
$\sigma^{(1)}\sigma^{(2)}=0$.  Without loss of generality,
we can assume that $\sigma^{(1)}=0$, otherwise, if $\sigma^{(2)}=0$  we can
use the linear transformation $x_1=y,\ y_1=x$,  which  replaces the points
$M_1$ and  $M_2$  and keeps the canonical form of system (35). Thus, after the
respective changing of parameters, we derive  the same system, but with
$\sigma^{(1)}=0$. In this case from (37) we obtain $c=eu+2m$ and for the
system
(35) we have
$$
\eqalign{
  & \mu= 4u(eh-eum-2m^2)^2,\cr
  & C_3={4\over 3}u(eh-eum-2m^2)(eh+3em+6emu+2e^2u^2+2e^2u+4m^2).\cr
}
$$
   Since $\mu\ne0$ from $C_3=0$ it results
$$
 e(h+3m+6mu+2eu^2+2eu)+4m^2=0.
$$
Therefore, we have $e\ne0$, otherwise it follows $m=0$ and, hence, $\mu=0$.
Thus, we can assume that $e=1$ (by scaling  time  if
necessary) and we get the relation: $h=-3m-6mu-2u^2-2u-4m^2$. In this
case after a shift of the origin of coordinates to the singular point
$M_1(1,0)$ of system (35) we get the system
$$
\eqalign{
  {dx\over dt}&= -(u+2m)x + [((u+2m)u+2\bar h]y - (u+2m)x^2 + 2\bar hxy -
                 (u+2m)uy^2,\cr
  {dy\over dt}&= -x + (u+2m)y - x^2 + 2mxy - uy^2
}\eqno(39)
$$
where $\bar h=-3m-6mu-2u^2-2u-4m^2$. From Proposition 3 applied to the
system (39) we obtain that
$$
\eqalign{
    I_1&=I_6=0,\ I_2 = 4(2u+3m)(u+2m+1) ,\cr
    I_3&=(u+m)(u+2m)(u+2m+1)(5u+12m+9) ,\cr
    I_{13}&=-4u(u+m)^3(2u+3m)(u+2m+1)^2,\cr
    5I_3-2I_4&= (u+m)(u+2m+1)(5u+12m+9)(7u+12m),\cr
   13I_3-10I_5&=-(u+2m+1)(48m^3+32m^2u+36m^2+9mu^2\cr
  & \quad      +9mu+5u^3-7u^2).
}\eqno(40)
$$
On the other hand  the following affine invariants can be
calculated for system (39)
$$
\eqalign{
  & \mu= 4u(2u+3m)^2(u+2m+1)^2,\cr
  &  C_4={1\over 3}I_{13},\  C_2= 3 I_2C_{4},\  C_5 = 6 I_3C_{4}.
}\eqno(41)
$$
It is easy to show,  from (40), (41) and $C_4\ne0$ that
conditions $I_2<0$ and $I_3=0$ are equivalent to $C_2C_4<0$ and $C_5=0$,
respectively. Note that conditions $5I_3-2I_4= 13I_3-10I_5=0$
can not be satisfied for system (39) because of $C_4\ne0$. Indeed, from
$5I_3-2I_4=0$ and (40) it follows $(5u+12m+9)(7u+12m)=0$, however, from (40)
it is easy to observe that in both determined by two factors cases we
have $13I_3-10I_5\ne0$.

 Notice also that given the indicated values of parameters $c, e$ and $h$
from $C_4\ne0$, (41) and (36) it follows that
$\sigma^{(2)}=-4(m+u)(2m+u+1)\ne0$.
 Hence, by [19] and Lemma 1  there exists
only one singular point of the center type for system (35).

 Thus, in the case $C_4\ne0$ the assertion of Theorem 3 is valid.
\bigskip
   2) Let us assume now $C_4=0$. In this case the condition $c+eu=0$ holds,
otherwise, as it was demonstrated above,  from $C_1=C_3=0$ we obtain that
$c=eu+2m,\ e=1,\ h=-3m-6mu-2u^2-2u-4m^2 $, and
$$
   \mu= 4u(2u+3m)^2(u+2m+1)^2,\
   C_4=-{4\over 3}u(u+m)^3(2u+3m)(u+2m+1)^2.
$$
Hereby in virtue of $\mu\ne0$ and $c+eu=2(m+u)\ne0$ it follows that
$C_4\ne0$.
This contradiction proves our assertion.

Thus, condition $c=-eu$ is satisfied and
hence, we arrive to the system
$$
\eqalign{
    {dx\over dt}&= -eux - eu^2y + eux^2 + 2hxy + eu^2y^2,\cr
    {dy\over dt}&=  ex + euy - ex^2 + 2mxy - euy^2,
}\eqno(42)
$$
for which
$$
   \mu = 4e^2u(h+mu)^2\ne0, \ C_1=0,\
   C_3=4e^2u(h+mu)(h-eu^2 -eu- mu).
$$
Therefore, from condition $C_3=0$ and $\mu\ne0$, we obtain
$h=eu^2 + eu+ mu$. In this case, for system (42) we also obtain
$$
\eqalign{
  \sigma^{(1)}& =2(eu +m),\ \mu  = 4e^2u^3(eu+e+2m)^2,\cr
  \sigma^{(2)}& =2u(eu+m),\
   C_8={16\over 3}e^2u^4(eu+m)^2(eu+e+2m)^2.
}\eqno(43)
$$
According to [19] for the existence of a center for system (42)
it is necessary that
$m=-eu$ and this condition together with  (43) and $\mu\ne0$,
is equivalent to $C_8=0$. In this case $\mu=4e^4u^3(u-1)^2$, and after shift
of the origin of coordinates to the singular point $M_1(1,0)$ system (42)
can be transformed into the system
$$
\eqalign{
    {dx\over dt}&=  eux + eu(2-u)y + eux^2 + 2euxy + eu^2y^2,\cr
    {dy\over dt}&=  -ex - euy - ex^2 - 2euxy - euy^2,
}
$$
for which
$$
   I_1^{(1)}= I_6^{(1)}=I_{13}^{(1)}=0 \ I_2^{(1)} = 4e^2u(u-1).
$$

   On the other hand, after a shift of the origin of coordinates to the
singular point $M_2(0,1)$ in the case under consideration system (42)
will be brought to the system
$$
\eqalign{
    {dx\over dt}&=  eux + eu^2y + eux^2 + euxy + eu^2y^2,\cr
    {dy\over dt}&=  e(1-2u)x -euy - ex^2 - 2euxy - euy^2,
}
$$
for which
$$
   I_1^{(2)}= I_6^{(2)}=I_{13}^{(2)}=0 \ I_2^{(2)} = 4e^2u^2(1-u).
$$
As it is easily seen $I_2^{(1)}I_2^{(2)}=-4\mu^2$. Therefore,
 if $\mu>0$ we get $I_2^{(1)}I_2^{(2)}<0$ and by Lemma 1 the
system (42) has only one singular point of center type.  If $\mu<0$ we obtain
$u<0$ and it is easy to see that in this case $I_2^{(1)}>0,\
I_2^{(2)}>0$. By  Lemma 1, system (42) has not any centers.

    Theorem 3 is proved.

\bigskip
{\bf Theorem 4.} {\sl For the existence of a center of system (1) with
$\mu\ne0, D=T=P=0,R\ne0$ (there are one simple and one triple singular
points) it is necessary and sufficient that the following conditions
hold:
$$
\eqalign{
   &  C_3=C_4=0,\ C_9>0;  \cr
   &  C_3=C_4=C_9=0,\ \mu>0.
}
$$
Moreover, the center is unique.}

\bigskip
{\bf Proof.} Let us assume that system (1) has one simple and one triple
singular points.  By applying the affine transformation we can move the
real  singular  points to the points $M_0(0,0)$ and $M_1(1,0)$, respectively.
Hence, system (1) becomes
$$
\eqalign{
    {dx\over dt}&= cx + dy - cx^2 + 2hxy + ky^2,\cr
    {dy\over dt}&= ex + fy - ex^2 + 2mxy + ny^2,
}\eqno(44)
$$
for which
$$
  \mu= (cn-ek)^2 + 4 (cm-eh)(hn-km)\ne0.
$$
Without loss of generality we can assume, that the critical point $M_0(0,0)$
is degenerated (of the third multiplicity). Therefore, the condition
$cf-de=0$ holds. For
$\mu\ne0$ it follows that
$c^2+e^2\ne0$ and one can sets $ c=d=0$. Indeed, these
relations can be obtained by the transformation $x_1=ex-cy,\ y_1=y$ if
$e\ne0$ and by the transformation $x_1=y,\ y_1=x$ if $e=0$ (in this case the
conditions $cf-de=0$ and $c\ne0$ imply $f=0$).

   Thus, for system (44) with $c=d=0$  the following comitants can be
calculated:
$$
     P=e^2(ek-2fh)^2(2hx+ky)^2y^2,\ \mu=e(ek^2-4h^2n+4hkm).
$$
From $\mu\ne0$ we obtain $e\ne0$ and we can assume $e=1$ (
scaling  the time  if  necessary) and, hence, the
condition $P=0$
implies $k=2fh$. Therefore, we get  the canonical form
$$
    {dx\over dt}= 2hxy + 2fhy^2,\
    {dy\over dt}= x + fy - x^2 + 2mxy + ny^2,\eqno(45)
$$
for which
$$
   \mu=4h^2(f^2+2fm-n),\  C_3=2hf(f^2+2fm-n). \eqno(46)
$$
According to Lemma 1 for the existence of a center it is necessary that
$C_3=0$.
From (46) and condition $\mu\ne 0$ we obtain $f=0$ and then
for this system $C_1=0$.
Hereby, after shifting of the origin of coordinates to the simple singular
point $M_1(1,0)$ of system (45) we get the following system
$$
    {dx\over dt}= 2hxy + 2fhy^2,\
    {dy\over dt}= x + fy - x^2 + 2mxy + ny^2,
\eqno(47)
$$
for which
$$
\eqalign{
   I_1&=m,\ I_2= 4(m^2-h),\ \mu  =-4h^2n,\cr
   C_4& ={2\over 3}mh(h+n)^2,\  C_9  =h[(h+n)^2 +m^2n].
}  \eqno(48)
$$
  If $C_4\ne0$ we obtain $I_1\ne0$ and according to Proposition 3, the
singular point (0,0) of system (47) (i.e. point $M_1(1,0)$ of system (45))
is not a center.

  Let us assume now that the
  condition $C_4=0$ is satisfied. By  (48)
and $\mu\ne0$ it follows that $m(h+n)=0$.
\smallskip
   1) If  $C_9>0$ the condition $m=0$ holds, otherwise for $h=-n$ we obtain
$C_9=-m^2n^2\le0$.  Therefore, $I_1=0$ and $I_2<0$ because of $Sgn I_2=
-Sgn h =-Sgn C_9$.

   Notice that for system (48) with $m=0$ the conditions
$I_6=I_{13}=0$ are satisfied and, hence, by Proposition 3 system (48)
has one center if $I_1=0,\ I_2<0$, i.e. if $C_4=0,\ C_9>0$.
\smallskip
   2) If  $C_9<0 $ by (48) and $\mu\ne0$ we obtain that either
$m=0,\ h<0$ (i.e. $I_1=0, \ I_2>0$) or $h=-n,\ -m^2n^2<0$ (i.e. $I_1\ne0$).
In both cases, by Proposition 3 system (47) has no  center.
\smallskip
   2) If  $C_9=0 $  the
   conditions $m(h+n)=0,\ (h+n)^2+m^2n=0$ and $hn\ne0$
imply  $m=0,\ n=-h$. Hereby, $\mu=4h^3$ and  $Sgn I_2=-Sgn \mu$. Therefore,
if $C_4=C_9=0,\ \mu>0$ for system (47) we obtain that
$I_1=I_6=I_{13}=0, \ I_2<0$
and hence, according to
Proposition 3, system (48) has a singular point of
a center type.

     Theorem 4 is proved.

\bigskip
\centerline{\bf \S 2. System with total multiplicity $m_f= 3$}
\bigskip

In this section we shall determine the conditions for the existence of a
center by using the
invariants (3) and Table 1 in the case where multiplicity
$m_f$ equals three.

From Table 1 and [19] it follows that the
system (1) with $m_f=3$ can has a center
only if it belong to set $M_{10}\cup M_{11}\cup M_{12}$. This implies 3 different
cases which will be examined in the sequel.

\bigskip

{\bf Theorem 5.} {\sl For the existence of a center of system (1) with
$\mu=0$, $H\ne0$, $D<0$ (there are three simple singular points)
it is necessary and sufficient that one of the
following two sequences of conditions holds:
$$
\eqalign{
   (i) \   & \ C_2C_4<0 ,\ C_1=C_3=C_5=0;\cr
   (ii) \, & \ C_4=0,\ C_1=C_3=C_{10}=0,\ C_{11}\le0.\cr
}
$$
Moreover, the center is unique.}
\bigskip

{\bf Proof.} Let us assume that system (1) has three simple
singular points.  By applying an affine transformation we can replace two
real  singular  points by the points $M_0(0,0)$ and $M_1(1,0)$, respectively.
Hence, system (1) will be transformed into the system
$$
\eqalign{
    {dx\over dt}&= cx + dy - cx^2 + 2hxy + ky^2,\cr
    {dy\over dt}&= ex + fy - ex^2 + 2mxy + ny^2,
}\eqno(49)
$$
for which, by using the notations
$$
\eqalign{
  \bar B=& cn-ek,\ \bar C= cm-eh\ \bar D=de-cf,\cr
  \bar E=& dn-fk,\ \bar F=fh-dm,\ \bar H=hn-km,
}
$$
 we obtain
$$
\eqalign{
      \mu &= \bar B^2 + 4 \bar C\bar H=0,\
       \tilde S  = -\bar C x^2-\bar Bxy +\bar Hy^2,\cr
       H &= [\bar C(\bar B +2\bar F)- \bar B\bar D]x+
         [\bar B(\bar B+\bar F)+ 2\bar H(\bar C+\bar D)]y\ne0.
}\eqno(50)
$$
We can see that if $\bar B=\bar C=\bar H=0$ then $H=0$ and the
conditions from
Theorem 5 are not valid. We shall prove, that in this case
for system (49) $\bar C\ne0$. Indeed, if we assume that
$\bar C=0$ from $\mu=0$ and (50) we
obtain $\bar B=0$ and, hence, $\bar H=(hn-km)\ne0$. Therefore, we get
the linear homogeneous system:
$$
 \bar C=cm-eh=0,\ \bar B=cn-ek=0,
$$
with determinant  (with respect to parameters $c$ and $e$) is equal
to $hn-km=\bar H\ne0$. Thus, $c=e=0$ and, hence, $\bar D=0$ and again $H=0$.
This proves our assertion.

As we can observe from (50), invariant $\mu$ is the discriminant of the
 comitant
$\tilde S$ and since $\mu=0$ we can write
$\tilde S=(\alpha x+\beta y)^2\ne0$, where
$\alpha\ne0$ for $\bar C\ne0$. Therefore, by applying the linear
transformation
$$
    x_1 = x+{\beta\over\alpha}y, \quad y_1=y,
$$
which keeps the
singular points $M_0$ and $M_1$, we obtain a new system of
the same form (49), for which comitant $\tilde S$ has the form
$\tilde S=\gamma x_1^2$
(see [15], p. 26). Taking into account (50) this form of comitant $K$ implies
the relations
$ ek-cn=hn-km=0$ and from $eh-cm\ne0$ we get  $n=k=0$. This leads
to the system
$$
    {dx\over dt}= cx + dy - cx^2 + 2hxy,
    {dy\over dt}= ex + fy - ex^2 + 2mxy ,\eqno(51)
$$
for which
$$
\eqalign{
   &\mu= 0,\  H=2\bar C \bar Fx\ne0, \cr
  & D=-{8\over 27}\bar D^2\bar F^2(2\bar C-\bar D)^2<0.\cr
 }\eqno(52)
$$
It is easy to see that system (51) has singular points
$$
    M_0(0,0),\quad M_1(1,0),\quad
    M_2\left({\bar D\over 2\bar C},
        {\bar D(2\bar C-\bar D)\over 4\bar C\bar F}\right)
$$
and from [19]   we obtain that  for these points
$$
\eqalign{
  \sigma^{(0)}&=c+f,\ \sigma^{(1)}=f+2m-c, \cr
  \sigma^{(2)}&= c+f + (m-c){\bar D\over \bar C} +
   h {\bar D(2\bar C-\bar D)\over 2\bar C\bar F},
} \eqno(53)
$$
correspondingly.
 For system (51) we have:
$$
     C_1=-24 \bar C\bar Fh\sigma^{(0)}\sigma^{(1)}\sigma^{(2)},
     \ C_4=  hP^\prime, \eqno(54)
$$
where $P^\prime$ is a polynomial in the coefficients of system (51),
and, hence,  if $C_1\ne0$  we get
$\sigma^{(0)}\sigma^{(1)}\sigma^{(2)}\ne0$ and by [19] system (51)
has no center.

    Let us assume that $C_1=0$. We
shall consider two cases: $C_4\ne0$ and $C_4=0$.
\medskip

   1) If $C_4\ne0$ then according to (54) $h\ne0$ and, hence, the condition
$C_1=0$ implies  that $\sigma^{(0)}\sigma^{(1)}\sigma^{(2)}=0$.
Without loss of generality we can assume that $\sigma^{(0)}=0$, otherwise
we can use a linear transformation, which will  replace the points $M_0$
and $M_1$  in  case $\sigma^{(1)}=0$ or points $M_0$ and $M_2$ in case
$\sigma^{(2)}=0$. In both cases, after the corresponding change of
parameters, we obtain the same system, but with $\sigma^{(0)}=0$. Thus,
from (53), we obtain $f=-c$ and for system  (51) the following invariants can
be calculated:
$$
\eqalign{
  & \bar D C_4=-{1\over 3}\bar C\bar Fh\sigma^{(1)}\sigma^{(2)}=
                                 {1\over 3}\bar DI_{13}^{(0)},\cr
  & C_2=3I_2^{(0)}C_4,\  C_3=-{2\over 3}I_{6}^{(0)},\ C_5=6I_3^{(0)}C_4,\cr
  & C_6=6(5I_3^{(0)}-2I_4^{(0)})C_4,\ C_7=2(13I_3^{(0)}-10I_5^{(0)})C_4,
}\eqno(55)
$$
where $I_j^{(0)}\ (j=2,3,4,5,6,13)$ are the values of the center affine
invariants of Proposition 3, calculated for system (51) with singular
point $M_0(0,0)$.  From $C_4\ne0$ and (54) we get
$\sigma^{(1)}\sigma^{(2)}\ne0$ and neither of the singular points
$M_1$ or $M_2$ can be a center. In this case, from
Proposition 3, relation $I_1=c+f=0$ and (55), we conclude that
the singular point
$M_0(0,0)$ will be a center if and only if the following conditions hold:
$$
     C_1=C_3=C_5(C_6^2+C_7^2)=0,\quad C_2C_4<0.
$$
However, we shall prove that in the case under consideration conditions
$C_5\ne0$ and $C_6=C_7=0$ can not  be satisfied.  Indeed, if we suppose
the contrary, then from (55) it follows:
$ I_{13}\ne0,\ I_2<0, \ I_1=I_6=0$, $I_3\ne0$ and $5I_3-2I_4=13I_3-10I_5=0$.
Hereby, as it was shown in [20, p.131], by applying a center affine
transformation  system (51) can be brought to the canonical system
$$
   {dx\over dt} =y+ qx^2+ xy,\quad {dy\over dt}=-x - x^2+3qxy + 2y^2,
$$
for which $D=8q^2(q^2+1)>0$ in virtue of $I_{13}=125q(q^2+1)/8\ne0$.
As we can observe this contradicts condition $D<0$ of Theorem 5.

\bigskip

  2) If the condition $C_4=0$ is satisfied from Lemma 1 (which are
  necessary for
the existence of a center)  conditions $C_1=C_3=0$ imply that
$h=0$. Indeed, if
$h\ne0$ as it was shown above, taking into account (54) and (55) we can
assume, without loss of generality, that conditions $C_1=C_4=0$ imply  that
$\sigma^{(0)}=\sigma^{(1)}=0$.
Hereby, from (53) it follows that $f=-c,\ m=c$ and for system (51) we
obtain $C_3={4\over 3}ch(c^2-eh)=0$, however
$H=-2c(c^2-eh)(h+d)x\ne0$.

   Thus,  $h=0$, and for system (51), using (52), we obtain
$$
\eqalign{
     H=2cdm^2x\ne0,\
    C_{10}= -cd^2m^3\sigma^{(0)}\sigma^{(1)}\sigma^{(2)}.
 }\eqno(56)
$$
 According to [19], for the existence of a center it is necessary that
 $C_{10}=0$
and from $H\ne0$ and (56), it follows that
$\sigma^{(0)}\sigma^{(1)}\sigma^{(2)}=0$. As it was mentioned above we can
assume $\sigma^{(0)}=0$. Therefore, the conditions $f=-c$ holds and for
system (51) we obtain
$$
\eqalign{
 &C_3={2\over 3}d(c-m)(c^2+2cm+de)=-{2\over 3}I_{6}^{(0)},\
   I_{1}^{(0)} = I_{13}^{(0)}=0,\cr
 &C_{11}= {8\over 3}d^2m^2(c-m)^2(c^2+de)=
   {4\over 3}d^2m^2(c-m)^2I_{2}^{(0)},
}\eqno(57)
$$
where $I_j^{(0)}\ (j=1,2,6,13)$ are the values of the center affine
invariants from Proposition 3, calculated for system (51) with
$h=0$ and $f=-c$.
\smallskip
  If $C_{11}\ne0$ we get from (57) that $Sgn C_{11}= Sgn I_{2}^{(0)}$
and from Proposition 3 and (57),  we conclude that the
singular point $M_0(0,0)$ of system (51) will be a center if and only
if $C_3=0$ and $C_{11}<0$.

   If condition $C_{11}=0$ is satisfied, then from $\bar D=de+c^2\ne0$
and (57), it follows that $m=c$, so our system becomes
$$
    {dx\over dt}= cx + dy - cx^2,
    {dy\over dt}= ex - cy - ex^2 + 2cxy ,\eqno(58)
$$
the singular points of which are the following
$$
   M_0(0,0),\quad M_1(1,0),\quad
   M_2\left({c^2+de\over 2c^2},\ {d^2e^2-c^4\over 4c^3d}\right).
$$
 According to Proposition 3, we shall calculate the values of the
 invariants $I_1,\ I_2,\ I_6$
and $I_{13}$ for system (58) as well as for others two   with critical
points $M_1$ and $M_2$ situated at the origin, respectively. Thus, we get
 the following expressions:
$$
\eqalign{
   & I_{1}^{(i)}=I_{6}^{(i)}=I_{13}^{(i)} =0 (i=0,1,2),\cr
   & I_{2}^{(0)}=2(c^2+de)\ne0,\ I_{2}^{(1)}=2(c^2-de)\ne0,\
   I_{2}^{(2)}={1\over 2c^2}(de+c^2)(de-c^2),
}\eqno(59)
$$
correspondingly.  Since the condition
$Sgn I_{2}^{(2)}= -Sgn(I_{2}^{(0)}I_{2}^{(1)})$ holds, we can conclude, that
among the quantities $I_{2}^{(i)}\ (i=0,1,2)$ one and only one will be
negative. Hence, in the case under consideration, one and only one of the singular
points of system (58) is of the center type.

  Since all possible cases were examined Theorem 5 is proved.
\bigskip

{\bf Theorem 6.} {\sl  For the existence of a center of system (1) with
$\mu=0$, $H\ne0$, $D>0$ (there are one simple real and two imaginary singular
points) it is necessary and sufficient that  the
following  conditions hold:
$$
    C_3=C_9=C_{10}=0, \  C_{11}<0.
$$
Moreover, the center is unique.}
\bigskip

{\bf Proof.} Let us assume that  system (1) has one simple real and two
imaginary singular points.  By applying the affine transformation we can
replace the real  singular  point by the origin.
Hence, system (1) becomes
$$
\eqalign{
    {dx\over dt}&= cx + dy + gx^2 + 2hxy + ky^2,\cr
    {dy\over dt}&= ex + fy + lx^2 + 2mxy + ny^2,
}\eqno(60)
$$
for which, by using the notations
$$
\eqalign{
 & \overset *\to B= gn-lk,\ \ \overset *\to C=cm-eh,\ \ \overset *\to D=de-cf,\cr
 & \overset *\to E= dn-fk,\ \ \overset *\to F=fh-dm,\ \ \overset *\to H=hn-km,\cr
 & \overset *\to G= gm-hl,\ \ \overset *\to L=dl-fg,\ \ \overset *\to N=cn-ek
}
$$
 we obtain
$$
\eqalign{
      \mu &=  \overset *\to B^2 - 4\overset *\to G\overset *\to H=0,\
       \tilde S  =  \overset *\to G x^2+\overset *\to Bxy +\overset *\to Hy^2,\cr
       H &=  [\overset *\to B(\overset *\to C -\overset *\to L)-
       2\overset *\to G(\overset *\to F+\overset *\to N)]x+
         [2\overset *\to H(\overset *\to C-\overset *\to L)-
         \overset *\to B(\overset *\to F+\overset *\to N)]y\ne0.
}\eqno(61)
$$
It is easy to see that from $\overset *\to H^2+\overset *\to G^2=0$ it
follows  that
$\overset *\to B=0$ and, hence, $H=0$, i.e. conditions of
Theorem 6 are not valid. Therefore, $\overset *\to H^2+\overset *\to G^2\ne0$
and we can consider $\overset *\to G\ne0$ (changing the coordinate
axes  if it is necessary).

As one can easily obtain from (61), the invariant $\mu$ is the discriminant
of the
comitant $\tilde S$ and in virtue of condition $\mu=0$ we can write
$\tilde S=(\alpha x+\beta y)^2\ne0$, where
$\alpha\ne0$ since $\bar C\ne0$. Therefore,  applying the linear
transformation
$$
    x_1 = x+{\beta\over\alpha}y, \quad y_1=y,
$$
we obtain a new system of
the same form (60), for which the comitant $\tilde S$ has the form
$\tilde S=\gamma x_1^2$
(see [15], p. 26]. By (61) this form of comitant $\tilde S$
implies the relations
$\overset *\to H=hn-km=0,\ \overset *\to B= gn-lk=0$ and from
$\overset *\to G=gm-hl\ne0$ we get  $n=k=0$. This leads
us to the system
$$
    {dx\over dt}= cx + dy + gx^2 + 2hxy,\
    {dy\over dt}= ex + fy + lx^2 + 2mxy ,\eqno(62)
$$
for which
$$
\eqalign{
 & H = 2\overset *\to G\overset *\to Fx\ne0,\
   D=-{8\over 27}\overset *\to D^2\overset *\to F^2 Z>0,\cr
 & Z=(\overset *\to L+2\overset *\to C)^2 + \overset *\to F\overset *\to M.
}\eqno(63)
$$
From $D>0$ it follows that $Z<0$ and
for system (62) we obtain
$$
\eqalign{
  & 2\overset *\to G\overset *\to F C_1  =
      3 h\sigma^{(0)}(R^2-ZI^2),\ \sigma^{(0)}=c+f,\cr
  & \overset *\to F C_4  ={h\over 12}[Rh+I(dG-2hC-3gF-4mF)],
    \ C_4= h^2\overset *\to G,
}\eqno(64)
$$
where
$$
\eqalign{
  R &= 4(c+f)\overset *\to F\overset *\to G +
      2(g+m)\overset *\to F(\overset *\to L-2\overset *\to C)+
        h(\overset *\to L^2+2\overset *\to F\overset *\to M+
        2\overset *\to D\overset *\to G),\cr
  I &= 2(g+m)\overset *\to F+h\overset *\to L.
}
$$
 There are two important cases to be examined: $C_9\ne0$ and $C_9=0$.
\bigskip
  I. Let us assume  first that condition $C_9\ne0$ hold. We shall
 prove the non-existence of a center for system (62).

     From $C_9\ne0$ and (64) we get $h\ne0$ and we can put
$h=1$ by changing  the time  if  necessary. Therefore,
we can consider, without loss of generality, that condition $g=0$ holds,
otherwise this can be obtained by applying the transformation
$x_1=x,\ y_1=gx/2+y$. Thus, we obtain  the system
$$
    {dx\over dt}= cx + dy +  2xy,\
    {dy\over dt}= ex + fy + lx^2 + 2mxy ,\eqno(65)
$$
for which in  accordance with Proposition 3, the condition $I_1=c+f=0$ gives
$f=-c$ and then
$$
     I_6=cdl+2cdm^2+2ce-d^2lm, \
   Z= 4c(cm^2-2cl-dlm-2em)+(dl-2e)^2<0.
$$
As condition $Z<0$ implies $c\ne0$ we shall introduce a new parameter
$u$ by setting $l=2cu$. Hereby, taking into consideration Proposition 3,
the necessary for the existence of the center  condition $I_6=0$
implies  that
$e=d^2mu-cdu-dm^2$ and for system (65) the following expressions can
be obtained:
$$
\eqalign{
  & Z=4(c-dm)[c(2du-m)^2-4c^2u + dm^2(3m-2du)]+4d^2m^2(du-2m)^2,\cr
  & H=-4cu(c+dm)x\ne0, \ I_1=I_6=0,\ I_{13}= 2u(c-dm)^2,\cr
  & 2I_2= (2c-d^2u)^2-d^2(2m-du)^2,\ I_3= 2c(2m-du), \cr
  & 5I_3-2I_4= 2(2m-du)(4c+dm),\cr
  & 13I_3-10I_5= 2(4c+dm)(3du+4m)+4dm(3m-4du).
}\eqno(66)
$$
\smallskip
  Since $H\ne0,\ Z<0$ and $I_2<0$ from (66) we get $I_3I_{13}\ne0$ and,
hence, from Proposition 3, for the existence of a center at the
origin of coordinates for system (65) it is necessary that conditions
$5I_3-2I_4=13I_3-10I_5=0$  be satisfied.
However, as it is easily seen
from (66) the
condition $5I_3-2I_4=0$ leads to the condition $c=-dm/4$, and,
hence,
$$
  13I_3-10I_5= 4dm(3m-4du)=0,\ 8I_2=5d^2m(4du-3m)<0.
$$
The obtained contradiction proves our assertion.
\bigskip

  II. Let conditions $C_9=0$ be satisfied, i.e. $h=0$. In this case for
system (62) we obtain $ H=2dgm^2x\ne0$ and, hence, the conditions $d=1$ and
$c=0$ can be considered  to be satisfied. Indeed, by $H\ne0$ it
results $d\ne0$ and by applying a change of the scale one can obtain
$d=1$. Hereby, after the
transformation $x_1=x,\ y_1=cx +y$ system (62) will be
transformed into the  system
$$
    {dx\over dt}= y + gx^2,\
    {dy\over dt}= ex + fy + lx^2 + 2mxy ,\eqno(67)
$$
for which
$$
\eqalign{
 & H = 2gm^2x\ne0,\  Z=(fg-l)^2+8egm<0,\cr
 & 3C_3=2[l(g+m)+fg(2m-g)],\cr
 & C_{10}=fm^2[f(lm-g^2f+gl)]-2e(g+m)^2],\cr
 & 3C_{11}=8m^2(g+m)[e(g+m)-lf]-4f^2gm^2(m-2g),\cr
 & I_1=f,\ I_2=2e+f^2,\ I_6=-(g+m)(l+fm),\ I_{13}=0.
}\eqno(68)
$$
  To prove the existence of a center in the case we discussing, according to
Proposition 3 it  remains  to prove the equivalence of the
sequence of the conditions  $C_3=C_{10}=0$, $C_{11}<0$ with the following one:
$I_1=I_6=0,\ I_2<0$.

    If conditions $C_3=C_{10}=0$ and $C_{11}<0$  are satisfied then $g+m\ne0$,
otherwise $m=-g$ and from (68) we obtain $C_{11}=4f^2g^4\ge0$.
Therefore some new parameter $w$ can be introduced,  namely:
$f= (g+m)w$. It is easily seen  from (68) that condition $C_3=0$
gives $l=gw(g-2m)$ and for system (67) we have
$$
  C_{10}=-2wm^2(g+m)^3(gmw^2+e)=0,\ Z=9g^2m^2w^2+8egm<0.
$$
We note that $e\ne -gmw^2$, otherwise $Z=g^2m^2w^2\ge0$. Hence, the condition
$C_{10}=0$ implies that $w=0$, i.e. $f=l=0$ and, according to (68),
from  $3C_{11}=8em^2(g+m)^2<0$ we get $e<0$.

    Thus, we have demonstrated, that the
    conditions $C_3=C_{10}=0$ and $C_{11}<0$
for system (67) with $Z<0$ are equivalent to $f=l=0,\ e<0$. On the other
hand, from (67) it is not difficult to observe, that
from $I_1=I_6=0,\ I_2<0$ and $Z<0$ it also follows that
$f=l=0,\ e<0$. This fact
proves our assertion.

  As all possible cases were examined, Theorem 6 is proved.
\bigskip

In accordance with [21] it occurs
\medskip

{\bf Theorem 7.} {\sl Quadratic system (1) with
$\mu=D=0$, $\ R\ne0$, $\ P\ne0$ (so there are one simple and one double
singular
points) has not a critical point of the center type.}
\bigskip

We shall include hear an independent {\bf proof} of this assertion.

 From Proposition 2 and
$\mu=D=0$ the conditions $R\ne0$ and $P\ne0$ are equivalent to $H\ne0$ and
$G^2-6HF\ne0$, respectively. We recall that in this case the
system (1) has one
simple and one double singular points
situated on its phase plane. In order to
find out the corresponding canonical form
of system (1) as it was indicated in the proof of Theorem 5,
by using an
affine transformation we can replace the two real singular points by
the points $M_0(0,0)$ and $M_1(1,0)$, respectively. Moreover, the
system (1) with conditions $\mu=0,\ H\ne0$ can be transformed by
 the linear transformation into the system
$$
    {dx\over dt}= cx + dy - cx^2 + 2hxy,\
    {dy\over dt}= ex + fy - ex^2 + 2mxy ,\eqno(69)
$$
for which, by using the notations
$$
   \bar C= cm-eh,\quad \bar D=de-cf,\quad \bar F=fh-dm,
$$
one can  obtain
$$
   \mu= 0,\ \ H=2\bar C \bar Fx\ne0,\ \
    D=-{8\over 27}\bar D^2\bar H^2(2\bar C-\bar D)^2.
$$
As it was shown above system (69) has singular points
$$
    M_0(0,0),\quad M_1(1,0),\quad
    M_2\left({\bar D\over 2\bar C},{\bar D(2\bar C-\bar D)\over 4\bar C\bar F}\right).
$$
Hereby, we can deduce that by $D=0$, the point $M_2$ will coincide
with $M_0$ for $\bar D=0$ and with point $M_1$ for $\bar D=2\bar C$.
     Without loss of generality we can assume that the condition $D=0$
implies $\bar D=de-cf=0$ and singular point $M_0$ becomes degenerated (the
points can be replaced if it is necessary). From
$\bar C\ne0$ we get $c^2+e^2\ne0$ and we can
set  $d=cu,\ f=eu$. Thus system (69) becomes the  system
$$
\eqalign{
    {dx\over dt}&= cx + cuy - cx^2 + 2hxy,\cr
    {dy\over dt}&= ex + euy - ex^2 + 2mxy,
}
$$
which, after placing of the origin at the point $M_1(1,0)$, gets into
the form
$$
\eqalign{
  &  {dx\over dt}= -cx + (cu+2h)y - cx^2 + 2hxy,\cr
  &  {dy\over dt}= -ex + (eu+2m)y - ex^2 + 2mxy.
}\eqno(70)
$$

According to Proposition 3 for the existence of a center it is necessary
$I_1=0$. Therefore, for system (70) condition $I_1= eu-c+2m=0$ yields
$c=eu+2m$ and then, for this system we obtain
$$
\eqalign{
   & I_6=2u(emu+2m^2-eh)(eh+2e^2u^2+6uem+4m^2)=0,\cr
   & I_2=  4(emu+2m^2-eh)<0.
}\eqno(71)
$$
From $I_2<0$, $I_6=0$ and (71) we get
$eh+2e^2u^2+6uem+4m^2=0$. Since $H\ne 0$, the  condition
$e^2+m^2\ne0$ holds. Then $e\ne0$ and,
by changing the time we can obtain $e=1$.
Therefore we have $h=-2(u^2+3um+2m^2)$ and for the system (70)
we obtain
$$
\eqalign{
   & I_1=I_6=0,\ I_2=4(2u+3m)(u+2m)\cr
   & I_3=(5u+12m)(u+m)(u+2m)^2,\cr
  & 5I_3-2I_4= (5u+12m)(7u+12m)(u+m)(u+2m),\cr
  & 13I_3-10I_5= -(u+2m)(5u^3+9mu^2+32m^2u+48m^3),\cr
  & H=-2u(2u+3m)^2(u+2m)^2x\ne0.\cr
}\eqno(72)
$$

It is easy to see that in virtue of $H\ne0$ and (72) the
condition $I_3=0$ yields $5u+12m=0$.  However, in this case
$I_2=72m^2/25\ge0$
and in accordance with Proposition 3 system (70) has not a center.

On the other hand, it is  easy to observe, that conditions
$5I_3-2I_4=13I_3-10I_5 =0$ can not be satisfied simultaneously,
because neither of the quantities $u=-12m/5$
nor $u=-12m/7$ will satisfy the relation $5u^3+9mu^2+32m^2u+48m^3=0$.
\smallskip

Thus, from Proposition 3, we can conclude that the affirmation of
Theorem 7 is valid.

\bigskip
\centerline{\bf \S 3. System with total multiplicity $m_f\le 2$}
\bigskip

In this section we shall determine the conditions for the existence of a
center in the case where the
total multiplicity $m_f$ is less than or equal to
two. But for a complete solution of the center problem, besides the
invariants (3) and those from Table 1 we also need the following
elements of the
minimal polynomial basis of the center-affine comitants [15]:
 $$
\eqalign{
 & I_{17}=   a^\alpha a^\beta_{\alpha \beta},\
 I_{20}=   a^\alpha a^\beta a^\gamma a^\delta_{\alpha \beta}
\varepsilon_{\delta \gamma},\cr
 &   K_1= a^\alpha_{\alpha\beta}x^\beta,\
    J_2= I_1(I_2-I_1^2)+4I_1I_{17}-4I_{20}.
}
 $$

\bigskip
{\bf Theorem 8.} {\sl System (1) with conditions $\mu=R=0,\ P\ne0,\
U>0$  (there are two simple singular points) has one center if and  only if
one of the following three sequences of conditions holds:
$$
\eqalign{
   (i)\   &\ C_2C_4<0,\ C_1=C_3=C_5=0;\cr
   (ii)\, &\ C_1=C_3=C_4=0,\ C_8>0;\cr
   (iii)  &\ \tilde S=0,\ K_1=0,\ I_1=0;\cr
}
$$
and it has two centers if and only if the following
conditions hold}
$$
\eqalign{
   (iv) &\ C_1=C_3=C_4=0,\ C_8<0,\ C_9>0.\cr
}
$$
\bigskip

{\bf Proof.} Let us assume that conditions  $\mu=R=0,\ P\ne0,
U>0$ are valid for system (1).  From Proposition 2,
we can easily  see that these conditions are eqiuvalent to the
following ones: $\mu=H=0,\ G\ne0, \ U>0$. Since in this case the
system (1) has
two real simple singular points, by using an affine
transformation we can replace them by the points $M_0(0,0)$ and $M_1(1,0)$,
respectively. Thus, we obtain the system
$$
\eqalign{
    {dx\over dt}&= cx + dy - cx^2 + 2hxy + ky^2,\cr
    {dy\over dt}&= ex + fy - ex^2 + 2mxy + ny^2,
}\eqno(73)
$$
for which, by using the notations
$$
\eqalign{
  \bar B=& cn-ek,\ \bar C= cm-eh\ \bar D=de-cf,\cr
  \bar E=& dn-fk,\ \bar F=fh-dm,\ \bar H=hn-km
}
$$
we obtain
$$
\eqalign{
      \mu &= \bar B^2 + 4 \bar C\bar H=0,\
       \tilde S = -\bar C x^2-\bar Bxy +\bar Hy^2,\cr
       H &= [\bar C(\bar B +2\bar F)- \bar B\bar D]x+
         [\bar B(\bar B+\bar F)+ 2\bar H(\bar C+\bar D)]y\ne0.
}\eqno(74)
$$
It will be convenient to examine the
cases: $\tilde S\ne0$ and $\tilde S=0$.
\bigskip

  {\bf A.} If $\tilde S\ne0$, it follows from the  proof of Theorem 5
and $\mu=0$, that
there exists a center affine transformation that brings
system (73) to the same form but with the
additional conditions: $k=n=0$. Thereby
for this system we obtain
$$
\mu=0,\ S= (cm-eh)x^2\ne0,\ H= 2(cm-eh)(fh-dm)x=0.
$$
Since $d^2+f^2\ne0$ (otherwise system (73) with $k=n=0$ becomes degenerated,
and then $G=0$) by  $H=0,\ \tilde S\ne0$, we may
assume, without los of generality, that conditions $h=du,\ m=fu$ hold. Thus,
in the case  $\tilde S\ne0$, we have obtained the following canonical form
$$
\eqalign{
    {dx\over dt}&= cx + dy - cx^2 + 2duxy,\cr
    {dy\over dt}&= ex + fy - ex^2 + 2fuxy,
}\eqno(75)
$$
for which
$$
\eqalign{
  & \tilde S=u(de-cf)x^2\ne0,\ G=(cf-de)^2(2u+1)x^2\ne0.\cr
}
$$
  On the other hand we shall simultaneously consider  the system
$$
\eqalign{
    {dx\over dt}&= -cx + d(2u+1)y - cx^2 + 2duxy,\cr
    {dy\over dt}&= -ex + f(2u+1)y - ex^2 + 2fuxy,
}\eqno(76)
$$
which is obtained from (75) by replacing the origin at the singular point
$M_1(1,0)$.

    For system (75), as well as for system (76), the values of the following
affine invariants can be calculated:
$$
\eqalign{
 & C_1= 12d^2u^2(2u+1)(cf-de)^2\sigma^{(0)}\sigma^{(1)},\cr
 & C_8= -{4\over 3}d^2u^2(2u+1)(de-cf)^2.\cr
}\eqno(77)
$$
We note, that the quantities $\sigma^{(0)}=c+f$ and $\sigma^{(1)}=-c+f+2fu$ correspond to
the singular points $M_0$ and $M_1$, respectively.
\smallskip

  I. Let us consider at first that condition $C_8\ne0$ is valid. According
to Lemma 1 from
$GC_8\ne0$ and (77) the condition $C_1 = 0$
(which is necessary for the existence of a center )
yields $\sigma^{(0)}\sigma^{(1)}=0$.
Without loss of generality one can consider $\sigma^{(0)}=0$ otherwise
the linear transformation  which replace the points $M_0$ and $M_1$ and
keeps the canonical form of system (75) can be applied.

   Thus, $f=-c$ and for system (75) we obtain
$$
\eqalign{
 & C_3=  {2\over 3}cd(2u+1)(de+c^2)(1-u),\cr
 & C_4= -{1\over 3}cd^2u^2(2u+1)(de+c^2)(u+1).\cr
}\eqno(78)
$$

   1) If $C_4\ne 0$ we shall prove that there exists no center for system
(75). Indeed, by (78) and Lemma 1 the second necessary for the existence of
a center condition $C_3=0$ by virtue of $C_4\ne0$ implies $u=1$. Therefore,
for system (75) in view of Proposition 3, we can  calculate
$$
\eqalign{
  & I_1^{(0)}=I_6^{(0)}=0,\
    I_{3}^{(0)}= d(de-8c^2)={1\over 3}(5I_3^{(0)}-2I_4^{(0)}),\cr
  & I_2^{(0)}=2(c^2+de),\
    I_{13}^{(0)}=-6cd^2(c^2+de).
}
$$
Since $I_{13}^{(0)}\ne0$ for $C_4\ne0$, according to Proposition 3, the
singular
point $M_0$ will be a center only if $I_{3}^{(0)}=0$, but the condition
$de=8c^2$ yields $I_2^{(0)}=18c^2>0$. Hence, from
$\sigma^{(1)}=-4c\ne0$ we deduce that for $C_4\ne0$ and $\tilde S\ne0$ the
system
(73) has no a center on its phase plane.
\smallskip

   2) Let us now assume that $C_4=0$. From $C_3=C_4=0$ and (78)
   we get $c=0$ and
for system (75) we obtain
$$
\eqalign{
  & I_1^{(0)}=I_6^{(0)}=I_{13}^{(0)}=0,\ I_2^{(0)}=2de,\cr
  & C_8={4\over 3}d^2e^2u^2(2u+1),\ C_9=d^3eu^3,
}
$$
and for system (76), in the same way we can  calculate
$$
   I_1^{(1)}=I_6^{(1)}=I_{13}^{(1)}=0,\ I_2^{(1)}=-2de(2u+1).
$$

Since the following relations
$$
 Sgn C_8=-Sgn\left(I_2^{(0)}I_2^{(1)}\right), \ Sgn C_9=-Sgn\left(I_2^{(0)}+I_2^{(1)}\right)
$$
hold,  from Proposition 3, we conclude that if $\tilde SC_8\ne0,$ then $
C_1=C_3=C_4=0$ and the
system (73) has one center if $C_8>0$, and it has two centers if
$C_8<0,\ C_9>0$.
\smallskip

  II. Let condition $C_8=0$ be satisfied. Taking into consideration
condition $\tilde SG\ne0$ and (77) condition $C_8=0$ yields $d=0$ and then, for
systems (75) and (76), we get  $ I_2^{(0)}= c^2+f^2\ge0$ and
$I_2^{(1)}=c^2+f^2(2u+1)^2$, respectively. By virtue of Proposition 3
there exists no center on the phase plane of system (73).
\bigskip

  {\bf B.} Let us assume that condition $\tilde S=0$ holds. From (74) it
  follows that
$ \bar C= cm-eh=0,\   \bar B= cn-ek=0,\  \bar H=hn-km=0$.
It is easy to see that condition $c^2+e^2\ne0$ holds for system (74),
otherwise this system becomes degenerated and then $G=0$. Therefore, without
loss of generality, one can be assumed that relations $h=cu,\ m=eu$ and
$k=cv,\ n=ev$ are valid, where $u$ and $v$ are some new independent
parameters.

    Thus, system (73) will be transformed into the system
$$
\eqalign{
    {dx\over dt}&= cx + dy - cx^2 + 2cuxy + cvy^2,\cr
    {dy\over dt}&= ex + fy - ex^2 + 2euxy + evy^2,
}\eqno(79)
$$
for which
$$
    \mu=0,\ \tilde S=0,\ G=-(cf-de)^2(-x^2+2uxy+vy^2)\ne0,\ K_1=(eu-c)x+(cu+ev)y.
$$

By translating the origin of coordinate at the singular point $M_1(1,0)$
of system (79) we obtain the system
$$
\eqalign{
    {dx\over dt}&= -cx + (2cu+d)y - cx^2 + 2cuxy + cvy^2,\cr
    {dy\over dt}&= -ex + (2eu+f)y - ex^2 + 2euxy + evy^2.
}\eqno(80)
$$
For system (79) as well as for system (80) one can find out the values of
the following affine invariants:
$$
\eqalign{
 & C_1= 12(cf-de)^2(v+u^2)(2ceu-c^2+e^2v)\sigma^{(0)}\sigma^{(1)},\cr
 & C_4= {1\over 3}(cf-de)(v+u^2)(2ceu-c^2+e^2v)(c-ue),\cr
 & C_8= {4\over 3}(cf-de)^2(v+u^2)(2ceu-c^2+e^2v),\cr
}\eqno(81)
$$
where $ \sigma^{(0)}=c+f$ and $ \sigma^{(1)}=-c+f+2eu$.
\smallskip

  I. If $C_8\ne0$ in accordance with Lemma 1, the
condition $C_1=0$  $G\ne0$ and (81) imply that
$\sigma^{(0)}\sigma^{(1)}=0$. By the same reasoning
we obtain $\sigma^{(0)}=0$. Hence $f=-c$, and for the
system (79) we obtain
$$
\eqalign{
 & C_3= -{1\over 3}(c^2+de)(2cu+d+ev)\sigma^{(1)},\cr
 & C_4= {1\over 6}(c^2+de)(v+u^2)(2ceu-c^2+e^2v)\sigma^{(1)},\cr
}\eqno(82)
$$
\smallskip
  1) Let us assume that the
  condition $C_4\ne0$ holds. Hence, $\sigma^{(1)}\ne0$
and singular point
$M_1(1,0)$ of system (79) (i.e. point (0,0) of system (80)) is not a center.
Hereby,   $C_3=0$ and (82) yield $d=-2cu-ev$ and we obtain:
$$
   C_2=-6(2ceu-c^2+e^2v)C_4,\ C_5=0.
$$
On the other hand, by Proposition 3 for system (79) we can  calculate
$$
   I_1^{(0)}=I_3^{(0)}=I_6^{(0)}=0,\ I_2^{(0)}=2(c^2-2ceu-e^2v)
$$
and, as it can be easily  seen, the following relation holds:
$Sgn(C_2C_4)=SgnI_2^{(0)}$. Hence, in accordance with Proposition 3 if
$C_2C_4<0$ there exists one center on the phase plane of system (79).
\smallskip

  2) If condition $C_4=0$ holds then since $C_8\ne0,
\ \sigma^{(0)}=0 $ and
(81) from (82) we get $\sigma^{(1)}= -c + f + 2eu =0$
and from $\sigma^{(0)} = c + f= 0$ we obtain $\sigma^{(1)}=2(eu-c)=0$
. Thus, $c=eu$ and for
system (79), as well as for system (80), we obtain
$$
\eqalign{
  & I_1^{(0)}=I_6^{(0)}=I_{13}^{(0)}=0,\ I_2^{(0)}=2(de+eu^2),\cr
  & I_1^{(1)}=I_6^{(1)}=I_{13}^{(1)}=0,\ I_2^{(1)}=-2(de+eu^2),\
}
$$
respectively. Since $I_2^{(0)}I_2^{(1)}<0$ by Proposition 3 for system
(73) either singular point $M_0$ or $M_1$ is of the center type.

  It is important to underline, that in the  case when $\tilde S=0,\ C_8\ne0$ and
$ C_4=0$
we have obtained $C_8={4\over 3}e^4(d+eu^2)^2(v+u^2)^2>0$ and, hence, these
condition can be united with the same one of the case $S\ne0$ excluding
conditions $\tilde S\ne0$ and $\tilde S=0$, and namely:

{\sl   If $C_8\ne0$ and $C_4=0$ system (73) has one center for
$C_1=C_3=0,\ C_8>0$ and has two centers for $C_1=C_3=0,\ C_8<0,\ C_9>0$}.

\bigskip
     II. Let us assume now that  condition $C_8=0$ hold, i.e. in accordance
with (81) one have $(v+u^2)(2ceu-c^2+e^2v)=0$. We shall prove that the
 condition $C_3=0$
 ( which is necessary for the existence of a center)
 implies $v+u^2=0$.
Indeed, if we assume  that $v+u^2\ne0$,  by (81), condition
$C_8=0$ yields $2ceu-c^2+e^2v=0$. We have already noted, that condition
$c^2+e^2\ne0$ is valid for system (80) and, hence, $e\ne0$ because of
$2ceu-c^2+e^2v=0$. Therefore, one can  consider $e=1$ (by changing
 the time  if it is necessary) and $v=c^2-2cu$. Hereby,
for system (79) we obtain
$$
C_3={2\over 3}(cf-d)^2(c-u),\ G=(cf-d)^2(x-cy)[x+(c-2u)y]\ne0.
$$
Hence, by  $G\ne0$, the condition $C_3=0$ yields $c=u$, and the relation
$v=c^2-2cu$ becomes  $v=-u^2$.

   Thus, conditions $C_8=C_3=0$ imply $v=-u^2$, and for system (79) we
can  calculate:
$$
   C_3={2\over 3}(cf-de)(c-ue)(fu-cu-d-eu^2),\ K_1=(eu-c)(x-uy).
$$
\smallskip

   1) If $K_1\ne0$, from $C_3=0$ it follows that
   $d=fu-cu-eu^2$ and for system
(79) as well as for system (80), we obtain
$$
    I_2^{(0)}=I_2^{(1)}=(c-eu)^2+ (f+eu)^2\ge0.
$$
Hence, according to Proposition 3 for $K_1\ne0$ there exists no center for the
system (79).
\smallskip

   1) If $K_1=0$ then $c=eu$ (this implies $C_3=0$) and for the
   systems (79) and
(80), respectively, we can obtain the following values of the invariants
$$
\eqalign{
  & I_1^{(0)}=f+eu,\ I_6^{(0)}=I_{13}^{(0)}=0,\
    I_2^{(0)}=2e(d-fu)+(f+eu)^2,\cr
  & I_1^{(1)}=f+eu,\ I_6^{(1)}=I_{13}^{(1)}=0,\
    I_2^{(0)}=-2e(d-fu)+(f+eu)^2.\cr
}
$$

As it can be easily  seen from $I_1^{(0)}=0$ we get
$I_1^{(1)}=0$ and
$I_2^{(0)}I_2^{(1)}<0$. Since for $K_1=0$ the
condition $I_1=0$ becomes an affine
invariant one, we can deduce, that in case of $K_1=0$ the
system (79) has one
singular point of the center type if and only if $I_1=0$.
\smallskip

 It remains to underline the following two moments.

   The $1^{st}$: For system (79) from $K_1=(eu-c)x+(cu+ev)y=0$ it follows that
$c=eu,\ v=-u^2$
and, hence, $C_3=C_8=0$, i.e. these last conditions can be excluded from
the respective consequence $\tilde S=0,\ C_8=C_3=0,\ K_1=0, I_1=0$.

   The $2^{nd}$:  If $\tilde S\ne0$  the conditions $C_2C_4<0,\
C_1=C_3=C_5=0$  are not compatible. Indeed, as it was indicated above (see
p. {\bf A}, I, 1)), if $C_4\ne0$ the
conditions $C_1=C_3=0$ yields $f=-c,\ u=1$ and
then for system (75) one can be get out:
$$
   C_4=-2cd^2(c^2+de)\ne0,\ C_2=-12cd^2(c^2+de)^2,\
   C_5= 6d(de-8c^2)C_4.
$$
Therefore, condition $C_5=0$ yields $de=c^2$ and, hence,
$C_2C_4=216c^8d^4>0$. This proves our assertion and we can exclude condition
$\tilde S=0$ from the indicated above sequence of conditions.

   As all cases were examined, Theorem 8 is proved.

\bigskip

{\bf Theorem 9.} {\sl For the existence of a center of system (1) with
$\mu=H=G=0,\ F\ne0$ (there are one simple real singular point) it is necessary and sufficient that one of the
following two sequences of conditions holds:}
$$
\eqalign{
   (i) &\ \tilde N\ne0,\ C_3=C_{10}=0,\ C_{11}<0;\cr
   (ii)&\ \tilde N=0,  \ J_2=0,\ J_1>0.\cr
}
$$
\bigskip

{\bf Proof.} Let us consider that conditions  $\mu=R=P=0$ and $U\ne0$
are valid for system (1).  Taking into consideration Proposition 2
one can easily be seen that these conditions are equivalent to the
following ones: $\mu=H=G=0,\ F\ne0$. Since in this case system (1) has
one real simple singular point we can consider it situated at the
origin of coordinates, i.e. system (1) will be of the form
$$
\eqalign{
    {dx\over dt}&= cx + dy + gx^2 + 2hxy + ky^2,\cr
    {dy\over dt}&= ex + fy + lx^2 + 2mxy + ny^2,
}\eqno(83)
$$
for which,
$$
\eqalign{
    &  \mu =  (gn-kl)^2-4(gm-hl)(hn-km)=0,\cr
    &   \tilde S = (gm-hl)x^2+(gn-kl)xy +(hn-km)y^2.
}\eqno(84)
$$
We shall  examine two cases: $\tilde S\ne0$ and
$\tilde S=0$.
\bigskip

  {\bf A.} If $\tilde S\ne0$  we will prove that
  point $M_0$ of system (83) is not
 of the center type. Indeed, the invariant $\mu$ is
the discriminant of the binary form $S$ and, hence, the
comitant $S$ can be represented
in the form $S=(\alpha x +\beta y)^2$. Since $S\ne0$ we can assume
$\alpha\ne0$, otherwise   replace $x$ and $y$.
Hereby, applying the linear transformation $x_1=\alpha x +\beta y,\ y_1=y$
we obtain  system (83) for which the following conditions hold:
$$
   gn-kl=hn-km=0,\  gm-hl\ne0.
$$
As it can be easily seen, these conditions yield $k=n=0$ and for system
(83)  we obtain:
$$
   S = (gm-hl)x^2\ne0,\ H=2(gm-hl)(dm-fh)x=0.
$$
Since $h^2+m^2\ne0$ for $S\ne0$ from $H=0$, without loss of generality,
one can put $d=2hu,\ f=2mu$, where $u$ is a new parameter (a constant factor
2 is introduced for computational considerations). This takes us
to the system
$$
\eqalign{
    {dx\over dt}&= cx + 2huy + gx^2 + 2hxy,\cr
    {dy\over dt}&= ex + 2muy + lx^2 + 2mxy,
}\eqno(85)
$$
for which,
$$
\eqalign{
  &    \mu = H=0,\quad S=(gm-hl)x^2\ne0,\cr
  &    G= 4u(gm-hl)[u(gm-hl)+(eh-cm)]x^2=0.
} \eqno(86)
$$

      1) If $h\ne0$ one can consider $h=1$ (otherwise a change of
scale can be done) and, hence, condition $G=0$ yields
$ e= cm-gmu+lu$. Hereby, from Proposition 3, for system (85)
we obtain $I_1=c+2mu=0$ and from this the following values of the
comitants can be
obtained:
$$
   F=2u^3(g+2m)(gm-l)^3x^3\ne0,\
   I_6=6u^2(g+2m)(gm-l).
$$
Therefore, condition $F\ne0$ yields $I_6\ne0$ and, from Proposition 3, the
singular point $M_0(0,0)$ of system (85) is not a center.
\smallskip

    2) If the condition $h=0$ is satisfied, from   $G=0,\ S\ne0$ and (80) it
follows that
$c=gu$. Hence, for system (85) we obtain $I_2=u^2(g^2+4m^2)\ge0$ and
again from Proposition 3, there is no center on the phase plane of
system (85).
\bigskip

  {\bf B.} Let us assume now that condition $\tilde S=0$ holds. From
  (84) we get
$gm-hl=gn-kl=hn-km=0$ and, hence, the homogeneous quadratic parts of system
(83) are proportional. Without any loss of generality we can assume
$g=h=k=0$, otherwise this can be obtained by applying a linear
transformation.

     Thus, system (83) can be transformed into the system
$$
    {dx\over dt}= cx + dy,\
    {dy\over dt}= ex + fy + lx^2 + 2mxy + ny^2,
\eqno(87)
$$
for which
$$
\eqalign{
  & \mu=H=0,\ G=(c^2n -2cdm+d^2l)(lx^2+2mxy+ny^2)=0,\cr
  & F=(cf-de)(cx+dy)(lx^2+2mxy+ny^2)\ne0,\ \tilde N= (m^2-ln)x^2.
} \eqno(88)
$$

   1) Let us assume
   that condition $\tilde N\ne0$ holds. Hereby, we shall prove that the
singular point $M_0(0,0)$ of the
system (87) will be of the center type if and
only if conditions $ C_3=C_{10}=0,\ C_{11}<0$ are satisfied.

   We will show  first that these conditions imply $d\ne0$, because one can
easily  see that, for $d=0$, system (87) has not a center at the origin.
Indeed, if we assume $d=0$ then, from (88) $G=0$ and $ F\ne0$
it results $n=0$ and for system (87) it follows at once $C_{11}=0$.

    Thus, $d\ne0$ and one can assume $d=1$, whence it follows from
$G=0$ and (88) that condition $l=2cm-c^2n$ holds. Hereby, for system (87)
we have
$$
\eqalign{
 & C_3=  {2\over 3}n(m-cn)(cf-e),\ \tilde N=(m-cn)^2x^2\ne0,\cr
 & F=(cf-e)(cx+y)(lx^2 + 2mxy + ny^2)\ne0.
}
$$
According to Lemma 1, in order to have a center on the phase plane for the
system (87) is
necessary that $C_3=0$, and by virtue
of $G\tilde N\ne0$ it follows at once that
$n=0$. Thus, system (87) becomes
$$
    {dx\over dt}= cx + y,\
    {dy\over dt}= ex + fy + 2cmx^2 + 2mxy,
\eqno(89)
$$
for which
$$
\eqalign{
  & C_{10}= 2m^4(c+f)(cf-e),\ C_{11}= {8\over 3}m^4(e-cf),\ I_1=(c+f),\cr
  & I_2=(c+f)^2+2(e-cf),\ I_6=-m^2(c+f),\  I_{13}=0.
} \eqno(90)
$$
Taking into consideration Proposition 3 and (90), we can deduce that  the
singular point $M_0(0,0)$ of system (89) will be a center if and only if the
conditions $I_1=0$ and $I_2<0$ are valid. By (90)
and $F\tilde N\ne0$ these conditions are equivalent to $C_{10}=0$ and
$C_{11}<0$, respectively. This  proves our assertion.
\bigskip

   2) Let us now assume that $\tilde N=0$, i.e. according to (88) the condition $m^2-ln=0$
is satisfied for system (87).
\smallskip
    a) If $d\ne0$ then by the same reason given above, we can put $d=1$ and
conditions $G=\tilde N=0$ from (88) yield $m=cn$ and $l=c^2n$. Thus, the
system (87)
becomes
$$
    {dx\over dt}= cx + y,\
    {dy\over dt}= ex + fy + c^2nx^2 + 2cnxy + ny^2,
\eqno(91)
$$
for which
$$
\eqalign{
  & \mu=H=G=0,\ F=n(cf-e)(cx+y)^3\ne0,\cr
  &  J_2=2(c+f)(e-cf),\ J_1= cf-e,\cr
  & I_1=c+f,\ I_2=(c+f)^2+2(e-cf),\ I_6=n^2(c+f)(e-cf).
} \eqno(92)
$$
By Proposition 3 and (92), we can conclude that the
singular point $M_0(0,0)$ of system (91) will be a center if and only if
conditions $I_1=0$ and $I_2<0$ are valid. Hereby, we can easily
see from (92) that  $F\ne0$ and these conditions are equivalent
to the following ones: $J_2=0,\ J_1>0$.

\smallskip
    b) In  case  $d=0$, it remain to be shown that conditions
$G=\tilde N=J_2=0$, and $J_1>0$ are not compatible. Indeed,
$d=0$ and $G=\tilde N=0$, imply $n=m=0$ and then it
follows at once that $J_1= cf, \ J_2=-2cf(c+f)$. Hereby it is not difficult to
see, that the condition $J_2=0$ yields $J_1\le0$ and this proves
our assertion.

     Theorem 9 is proved.

\bigskip
\centerline{{\bf REFERENCES}}

    [1] Dulac H. {\sl D$\acute e$\!t$\acute e$\!rmination et
int$\acute e$\!grations d'une cer\-taine classe d'$\acute e$\!quations
dif\-f$\acute e$\!ren\-tielles ayant pour point singulier un centre.}
Bull. sci. Math., 32 (1908), 230--252.
\smallskip

    [2] Kapteyn.W. {\sl On the centra of the integral curves which satisfy
differential equations of thq first order and the first degree.}
Proc. Kop. Akad. Wet., Amsterdam, {\bf 13}, 2(1911),1241--1252.
\smallskip

    [3] Kapteyn.W. {\sl New resheasha upon the centra of the integrals
which satisfy differential equations of the first order and the first degree.}
Proc. Kop. Akad. Wet., Amsterdam, {\bf 14}, 2(1912),1185--1195;\
{\bf 15}, 2(1912),46--52.
\smallskip

    [4] Frommer M. {\sl $\ddot U$\!ber das Auftreten von Wirbeln und
Strudeln in der Umgebung retionaler Unbestimmtheitsstellen.} Math. Ann.,
{\bf 109}, 3(1934), 395--424.
\smallskip

    [5] Saharnikov N.A.{\sl About Frommer's conditions for the existence of
a center.} Prikl. Mat. i Meh. {\bf 12}, 5(1948), 669--670 (Russian).
\smallskip

    [6] Sibirsky K.S.{\sl About conditions for the existence of
a center and a focus.} Uch. zap. Kishin. un-ta, 11(1954), 115-117 (Russian).
\smallskip

    [7] Belyustina L.N.{\sl About conditions for the existence of
a center.} Prikl.Mat. i Meh. {\bf 18}, 4(1954), 511 (Russian).
\smallskip

    [8] Malkin K.E.{\sl About some method on the center problem}.
Uch. zap. Ryazan. ped. in-ta, 24(1960), 107-117 (Russian).
\smallskip

   [9] Kukles I.S.{\sl Some criterions for the existence of a center.}
Tr. Samarc. un-ta, 47(1951), 29--98 (Russian).
\smallskip

    [10] Sibirsky K.S.{\sl The algebraic invariants of differential
equations and matrices}. "Shtiintsa", Kishinev, 1976.
\smallskip

   [11]  D.Bularas [Boularas], N.I.Vulpe, K.S.Sibirskii. {\sl The problem
of a center "in the large" for a general quadratic system}. Dokl. Akad.Hauk
SSSR, Tom 311 (1990), no. 4. Transl. Soviet Math. Dokl. Vol. 41 (1990), no. 2,
287--290.
\smallskip

   [12]  D.Bularas [Boularas], N.I.Vulpe, K.S.Sibirskii. {\sl The solution of
the  problem of a center "in the large" for a general quadratic differential
system}. Diff. uravnenia, V. 25(1989), no.11, 1856-1862 (Russian).
\smallskip

  [13]  V.A.Baltag and N.I. Vulpe. {\sl Number and multiplicity of singular
   points of a  quadratic differential system}. Ross. Akad. Nauk Dokl.,
 Tom 323 (1992), no. 1. Transl. Russian Acad. Sci. Dokl. Math. Vol. 45 (1992),
 no. 2 235--238.
\smallskip

 [14]  Baltag V.A., Vulpe N.I. {\sl Affine-invariant conditions for
 determining the number and multiplicity of singular points of
 quadratic differential systems }. Buletinul A.\c S. a R.M. ser. Matematica,
 1993, Vol. 1(11), 39--48 (English).
\smallskip

 [15]  Sibirsky K.S. {\sl Introduction to the Algebraic Theory of Invariants of
 Differential Equations}. Manchester University Press. Manchester, 1988.
\smallskip

 [16]  Boularas D., Calin Iu., Timochouk L., Vulpe N. {\sl
 $T$- comitants of quadratic systems: a study via the translation invariant}.
 Report of Faculty of Technical Mathematics and Informatics,
 Delft University of Technology, Report 96-90, Delft 1996.
\smallskip

 [17]  G.Gurevich, {\sl Foundation of the Theory of Algebraic Invariants}.
 Noordhoff, Groningen, 1964.
\smallskip

  [18]  Reyn J.W. {\sl Phase portraits of quadratic systems with finite
multiplicity one.} Nonlinear Analysis, Theory, Methods and Applicatins,
(to appear).
\smallskip

   [19]  Andronov A.A, Leontovich E.A., Gordon I.I. and Maier A.L.
{\sl Qualitative theory of second-order dynamical systems.}
    John Wiley \b Sons,  New York, 1973
\smallskip

   [20]  Vulpe N.I.{\sl Polynomial bases of comitants of differential
systems and their applications in qualitative theory},
"Shtiintsa", Chishinau, 1986 (Russian).
\smallskip

   [21] Berlinskii A.N.{\sl On the coexistence of singular points of
different types}. Izv. Vyss. Uchebn. Zaved. Mat. 1960, 5(18), pp. 27--32
(Russian).

\bigskip
\bigskip
{\it Institute of Matematics,}    \hfill      {\it 443 Eucaliptus Drive,}

{\it Academy of Science of Moldova}  \hfill  {\it Redlands, CA,}

{\it 5 Academiei Str, Chishin\u au,}  \hfill {\it 92373, U.S.A.}

{\it MD-20028, Moldova}  \hfill    {\it E-mail:}

{\it Phone:\ \ (373-2) 727059}  \hfill  {\it avoldman\@zzyzx.math.csusb.edu}

{\it Fax:\hskip 0.65cm (373-2) 738027}  \hfill

{\it E-mail: 15vulpe\@mathem.moldova.su}  \hfill

\enddocument